\documentclass[12pt]{amsart}
\usepackage{amsmath,amssymb,mathrsfs,amsfonts,mathtools}
\usepackage{graphicx,geometry}
\usepackage{dynkin-diagrams}
\usetikzlibrary{decorations.pathreplacing,fit,positioning}

\usepackage{setspace}

\usepackage[page wise]{ line no } 
\usepackage{physics}                    
\usepackage{titleps}

\usepackage{tikz}
\usepackage{tikz-cd}
\usepackage{amssymb}
\usepackage{verbatim}
\usepackage{textcomp}

\usepackage{float,comment,xcolor}
\usepackage{enumitem}
\usepackage{tikz-cd}
\usepackage{fancyhdr}
\newcommand{\alphalabel}[1]{\alpha_{#1}}

\graphicspath{{plot/}}
\newtheorem{theorem}{Theorem}[section]
\newtheorem{lemma}[theorem]{Lemma}
\newtheorem{example}[theorem]{Example}

\newtheorem{prop}[theorem]{Proposition}
\newtheorem{corollary}[theorem]{Corollary}

\theoremstyle{definition}
\newtheorem{definition}[theorem]{Definition}
\theoremstyle{remark}
\newtheorem{remark}[theorem]{Remark}
\numberwithin{equation}{section}
\def\nd{\noindent}

\begin{document}

\title[]{Counting rational curves and standard complex structures on hyperK\"ahler ALE 4-manifolds}

\author{Yuanjiu Lyu$^\dagger$}
\address{School of Mathematical  Sciences, University of Science and Technology of China, Hefei 230026 China}
\email{lyj999@mail.ustc.edu.cn}

\author{Bin Xu}
\address{CAS Wu Wen-tsun Key Laboratory of Mathematics and School of Mathematical  Sciences, University of Science and Technology of China, Hefei 230026 China}
\email{bxu@ustc.edu.cn}

\thanks{B.X. is supported in part by the Project of Stable Support for Youth Team in Basic Research Field, CAS (Grant No. YSBR-001) and NSFC (Grant Nos. 12271495 and 12071449). }

\thanks{$^\dagger$Y.L. is the corresponding author.}

\fancypagestyle{plain}{
  \fancyhead[LE,RO]{\shorttitle}
  \fancyfoot[C]{\thepage}
}

\date{\today}

\maketitle
\pagestyle{plain}

\begin{abstract}

\footnotesize All hyperK\"ahler ALE 4-manifolds with a given non-trivial finite group $1\not=\Gamma< {\rm SU}(2)$ at infinity are parameterized by an open dense subset of a real linear space of dimension $3{\rm rank} \,\Phi$. Here, $\Phi$ denotes the root system associated with $\Gamma$ via the McKay correspondence. Such manifolds are diffeomorphic to the minimal resolution of a Kleinian singularity. By using the period map of the twistor space,  we specify those points in the parameter space at which the hyperK\"ahlerian family of complex structures includes the complex structure of the minimal resolution. Furthermore, we count the rational curves lying on each hyperK\"ahler ALE 4-manifold. As an application, we prove that the twistor space of any hyperK\"ahler ALE cannot be K\"ahlerian. In particular, we strengthen some results of Kronheimer \scriptsize{{(J. Differential Geom., 29(3):665--683, 1989)}} \footnotesize and provide examples of non-compact and non-K\"ahlerian twistor spaces.

\vskip 4.5mm

\nd \begin{tabular}{@{}l@{ }p{10.1cm}} {\bf Keywords } &

hyperK\"ahler ALE,  Kleinian singularity, twistor space, root system,  rational curve.
\end{tabular}

\nd {\bf 2000 MR Subject Classification }
53C28, 17B22 

\end{abstract}


\section{Introduction}

\subsection{Background}
A \emph{hyperK\"ahler structure} on a Riemannian manifold $\left(M,g\right)$ is, by definition, a triple $(I,J,K)$ of three parallel complex structures on $M$ that  satisfies the quaternion relation $IJ = -JI = K$. The set $\left\{aI + bJ + cK\,|\,(a,b,c)\in\mathbb{S}^{2}\right\}$ is called \emph{the hyperK\"ahlerian family of complex structures} on $(M,g)$.
We say a 4-dimensional Riemannian manifold $(M,g)$ is \emph{asymptotically locally Euclidean (ALE)} if there exist a compact set $K \subset M$, a finite subgroup $1\not=\Gamma < \mathrm{SO}(4)$, and a diffeomorphism $\Psi: M \backslash K \to (\mathbb{R}^{4} \backslash B_{R})/\Gamma$, where $B_{R}$ is the standard ball of radius $R$ in $\mathbb{R}^4$, such that for each non-negative integer $k$, there holds
$$
|\nabla^{k}(\Psi_{*}g - g_{\rm Euc})(x)| = O(|x|^{-4-k})\,,\qquad x\in (\mathbb{R}^{4}\backslash B_{R})/\Gamma\qquad {\rm and}\qquad |x|\to +\infty.
$$
Kronheimer \cite{Kron1, Kron2} classified all hyperK\" ahler ALE 4-manifolds. He constructed explicitly for each finite group $1\not=\Gamma<\mathrm{SU}(2)$ a complete family of hyperK\"ahler ALEs with group $\Gamma$ at infinity. We provide a brief overview as follows. Let's fix a non-trivial finite group $\Gamma<\mathrm{SU}(2)$.

Denote $Q$ as the natural two-dimensional complex representation of  $\Gamma$ and $R$ as the regular representation. Define a linear space $$M(\Gamma) := \mathrm{Hom}_{\Gamma}\bigl(R,Q\otimes R\bigr)$$ equipped with a natural hyperK\"ahler structure when we identify $Q = \mathbb{C}^{2} \simeq \mathbb{H}$.
This space $M(\Gamma)$ is acted upon by a Lie group $\mathrm{U}(\Gamma)$, which preserves the hyperK\"ahler structure on $M(\Gamma)$. Denote by $Z$ the dual of the center of the Lie algebra $\mathfrak{U}(\Gamma)$ of $\mathrm{U}(\Gamma)$. There exists a natural isomorphism 
\begin{equation*}
    \rho : Z \xrightarrow{\sim} \mathfrak{h},
\end{equation*}
where $\mathfrak{h}$ is the Cartan subalgebra associated with $\Gamma$, see subsection 2.1 for $\mathfrak{h}$.
With respect to the $\mathrm{U}(\Gamma)$-action, we can associate a hyperK\"ahler moment map 
$$\mu_{\mathbb{H}}: M(\Gamma) \to \mathfrak{U}(\Gamma)^{*} \otimes \mathbb{R}^3$$ 
and consider the hyperK\"ahler quotient $X_{\zeta} := \mu_{\mathbb{H}}^{-1}(\zeta) / \mathrm{U}(\Gamma)$ for each $\zeta\in Z\otimes\mathbb{R}^{3}$. For details on hyperK\"ahler quotient, see \cite{HKLM}. 
These quotient manifolds may not be smooth. However, for $\zeta\in (Z \otimes \mathbb{R}^3)^\circ$,  $X_{\zeta}$ is a smooth ALE 4-manifold equipped with the quotient metric $g_{\zeta}$, where
\begin{equation*}
    (Z \otimes \mathbb{R}^3)^\circ  := (Z \otimes \mathbb{R}^3)\backslash \bigcup_{\theta\in\Phi}(D_{\theta}\otimes \mathbb{R}^3)\,,\qquad D_{\theta} := \mathrm{Ker}(\theta\circ\rho) \subset Z
\end{equation*}
and the union is taken over all roots in the root system $\Phi$ associated with $\mathfrak{h}$. By the property of hyperK\"ahler quotient, the smooth quotient Riemannian manifold $\left(X_{\zeta},g_{\zeta}\right)$ admits a natural hyperK\"ahler structure. We denote by $I_{\zeta}, J_{\zeta}, K_{\zeta}$ the three compatible complex structures on $X_{\zeta}$, and by $\omega_{\zeta}^k$, $k \in \{1, 2, 3\}$ the corresponding K\"ahler forms. The ordering of $I_{\zeta}, J_{\zeta}, K_{\zeta}$ is inherited from the ordering of complex structures on the hyperK\"ahler manifold $M(\Gamma)$. 

Furthermore, Kronheimer showed that each smooth manifold $X_{\zeta}$ for $\zeta\in(Z \otimes \mathbb{R}^3)^\circ$ is diffeomorphic to the minimal resolution of $\mathbb{C}^{2}/\Gamma$. We denote the minimal resolution of the $\mathbb{C}^{2}/\Gamma$ by $\eta : \widetilde{\mathbb{C}^2 / \Gamma} \to \mathbb{C}^2 /\Gamma$ (see Definition \ref{def_resolution_KLeinian}). The set $\eta^{-1}(0)$ is called the \emph{exceptional locus}, it is a union of $\mathrm{rank}\,\Phi$ rational curves.
We call the complex structure on the minimal resolution $\widetilde{\mathbb{C}^2 / \Gamma}$ the \emph{standard complex structure}, which we denote by $I_{\mathrm{std}}$.

\subsection{Main results}
For $\zeta = (\zeta^1,0,0) \in (Z \otimes \mathbb{R}^{3})^\circ$, Kronheimer \cite[Proposition 3.10]{Kron1} showed that the complex manifold $X_{\zeta}$, with respect to the complex structures $I_{\zeta}$, is isomorphic to $\widetilde{\mathbb{C}^2 / \Gamma}$. Our first main result characterizes all such $\zeta$  as follows. Denote by $\mathrm{rank}(\zeta)$ the dimension of  $\mathrm{span}_{\mathbb{R}}\left\{\zeta^{1},\zeta^{2},\zeta^{3}\right\}$. It is clear that $1\leq \mathrm{rank}(\zeta)\leq 3$ for each $\zeta\in (Z\otimes \mathbb{R}^3)^\circ$.

\begin{theorem}
\label{tm_1}
Let $\zeta\in(Z \otimes \mathbb{R}^{3})^\circ$. The hyperK\"ahlerian family of complex structures on $(X_\zeta, g_\zeta)$ includes
$I_{\rm std}$ if and only if ${\rm rank}(\zeta)=1$, that is $\zeta$ is of the form $\left(x \zeta^1, y \zeta^1, z \zeta^1\right)$ for some non-zero vector $(x,y,z)\in\mathbb{R}^{3}\setminus\{0\}$ and $\zeta^{1}\in Z\setminus\bigcup_{\theta\in\Phi}D_{\theta}$. Under this context, the hyperK\"ahlerian family of $X_{\zeta}$ includes exactly two complex structures which are equivalent to $I_{\rm std}$, and all the other complex structures are mutually equivalent and isomorphic to a hypersurface in  $\mathbb{C}^3$.
\end{theorem}

The \emph{twistor space} \cite[section 3.F]{HKLM} of the hyperK\"ahler manifold $\left(X_{\zeta},g_{\zeta}\right)$ is a complex threefold, denoted by $\mathrm{Z}(X_{\zeta})$, equipped with a holomorphic projection $\pi=\pi_\zeta: \mathrm{Z}(X_{\zeta}) \to \mathbb{P}^1$. For $u\in\mathbb{P}^{1} = \mathbb{C}\cup\{\infty\}$ we define a complex structure $I_{\zeta}^{u} := aI_{\zeta} + bJ_{\zeta} + cK_{\zeta}$, where
\begin{equation}
\label{eq_iden_S_P}
    (a,b,c) := \left(\frac{1-\norm{u}^2}{1+\norm{u}^2}, \frac{2\Re\,u}{1+\norm{u}^2}, \frac{-2\Im\,u}{1+\norm{u}^2}\right).
\end{equation}
Note that equation (\ref{eq_iden_S_P}) provides an identification between $\mathbb{P}^{1}=\mathbb{C}\cup\{\infty\}$ and 
$$\mathbb{S}^{2}=\left\{(a,b,c)\in \mathbb{R}^3|a^2+b^2+c^2=1\right\}.$$  
The fiber $\pi^{-1}(u)$ over $u$ is a complex submanifold of $\mathrm{Z}(X_{\zeta})$ that is isomorphic to the complex manifold $X_{\zeta}$ with respect to the complex structure $I_{\zeta}^{u}$. 
Our next result is:

\begin{theorem}
\label{tm_2}
 Fix a non-trivial finite subgroup $\Gamma$ in $\mathrm{SU(2)}$. Recall the linear space $Z$ and root system $\Phi$ associated with $\Gamma$. Given $\zeta\in(Z\otimes\mathbb{R}^{3})^\circ$, we could count rational curves fiberwise lying in the twistor space $\pi: \mathrm{Z}(X_{\zeta}) \to \mathbb{P}^{1}$ as follows{\rm :}
    \begin{itemize}
        \item[(i)] If $\mathrm{rank(\zeta)} = 1$, then there hold
            \begin{equation}
                \label{eq_1_rk1}
                \#\left\{u \in \mathbb{P}^{1} \,|\, \pi^{-1}(u) \mathrm{\,\,contains\,\, a\,\, rational\,\, curve} \right\} = 2,
            \end{equation}
and
            \begin{equation}
                \label{eq_1_rk1}
                \sum_{u\in\mathbb{P}^{1}}\#\left\{ \mathrm{rational\,\, curves\,\, in\,\,} \pi^{-1}(u) \right\} = 2r,
            \end{equation}
where $r := \mathrm{rank}\,\Phi$. 
        \item[(ii)] If $\mathrm{rank(\zeta)} = 2$, then there hold 
            \begin{equation}
                \label{eq_1_rk2}
                4 \leq \#\left\{u \in \mathbb{P}^{1} \,|\, \pi^{-1}(u) \mathrm{\,\,contains\,\, a\,\, rational\,\, curve} \right\} \leq \#\Phi, 
            \end{equation}
       and
            \begin{equation}
                \label{eq_2_rk2}
                2r\leq \sum_{u\in\mathbb{P}^{1}}\#\left\{ \mathrm{rational\,\, curves\,\, in\,\,} \pi^{-1}(u) \right\} \leq\# \Phi.
            \end{equation}
        \item[(iii)] If $\mathrm{rank(\zeta)} = 3$, then there hold 
            \begin{equation}
                \label{eq_1_rk3}
                3 \leq \#\left\{u \in \mathbb{P}^{1} \,|\, \pi^{-1}(u) \mathrm{\,\,contains\,\, a\,\, rational\,\, curve} \right\} \leq \#\Phi, 
            \end{equation}
       and
            \begin{equation}
                \label{eq_2_rk3}
                2r-1\leq \sum_{u\in\mathbb{P}^{1}}\#\left\{ \mathrm{rational\,\, curves\,\, in\,\,} \pi^{-1}(u) \right\} \leq\# \Phi.
            \end{equation}
    \end{itemize}
\end{theorem}

\begin{remark}
    The upper bounds for (\ref{eq_1_rk2}) and (\ref{eq_1_rk3}) were established by Kronheimer \cite[Lemma 2.7]{Kron2}. All lower bounds presented here are new results. However, determining optimal lower bounds remains challenging, primarily due to the underlying combinatorial problem discussed in Section 4. We plan to address this optimization problem in subsequent work.
\end{remark}

Considering all $\zeta\in\left(Z\otimes\mathbb{R}^3\right)^\circ$, we define a map as: 
\begin{align*}
    Q_1 : \left(Z\otimes\mathbb{R}^3\right)^\circ &\to \mathbb{Z}\\
    \zeta &\mapsto \#\left\{u \in \mathbb{P}^{1} \,|\, \pi_\zeta^{-1}(u) \mathrm{\,\,contains\,\, a\,\, rational\,\, curve} \right\}
\end{align*}
Although an explicit expression for $Q_1$ appears difficult to obtain, we establish:
\begin{theorem}
\label{tm_3}
    $Q_1$ is lower semi-continuous, that is, for every $\zeta_0\in\left(Z\otimes\mathbb{R}^3\right)^\circ$, there holds
    \begin{equation*}
        \liminf_{\zeta \to \zeta_0} Q_1(\zeta) \geq Q_1(\zeta_0)
    \end{equation*}
\end{theorem}

\begin{remark}
    Verifying the semi-continuity property for rational curve counts remains challenging.
\end{remark}

Combining Theorem \ref{tm_1} and Theorem \ref{tm_2}, we extend certain results of Hitchin \cite{Hit81} concerning the non-compact non-K\"ahlerian twistor spaces: 
\begin{corollary}
\label{cor_3}
For any finite group $1\not=\Gamma<\mathrm{SU(2)}$ and any $\zeta \in (Z \otimes \mathbb{R}^3)^\circ$, the twistor space $\mathrm{Z}(X_{\zeta})$ can not be K\"ahlerian (that is non-K\"ahlerian).
\end{corollary}

\subsection{Outline}
Section 2 reviews fundamental aspects of Kleinian singularities and their semi-universal deformations. The core technical tool, developed in Subsections 2.3 and 2.4, is the period map associated with the twistor space of hyperK\"ahler ALE manifolds.
Section 3 provides proofs for Theorems \ref{tm_1} and \ref{tm_2} in Case (i).
Section 4 completes the proof of Theorem \ref{tm_2} by reformulating it as a decomposition problem for irreducible root systems. This formulation naturally leads to new questions regarding root systems in Lie theory.
Section 5 establishes the lower semi-continuity of $Q_1$ and presents the proof of Theorem 1.4.
The final section derives Corollary \ref{cor_3} from Hitchin's foundational work \cite{Hit81}.

\subsection{Some basic notations}

\begin{itemize}
    \item $\left(X_{\zeta},I_{\zeta}^{u}\right)$: The complex manifold with underlying differential manifold $X_{\zeta}$ and complex structure $I_{\zeta}^{u}$. 
    \item $\mathfrak{h}^{\mathbb{C}}$: The complexification of the real Cartan subalgebra $\mathfrak{h}$.
    \item $\mathbb{P}(V)$ or $\mathbb{P}V$: The projective space associated with the complex vector space $V$.
    \item $\# P$: The cardinality of the set $P$.
    \item $P_{1}\sqcup P_{2}$: The disjoint union of two sets $P_{1}$ and $P_{2}$.
    \item $\mathrm{i}$: The imaginary unit, $\sqrt{-1}$.

\end{itemize}

\section{The Kleinian singularities}
Building on Kronheimer's foundational work \cite{Kron1,Kron2} establishing key connections between hyperK\"ahler ALE spaces and Kleinian singularities, this section reviews essential background and prepares for proving our main theorems. Specifically: Let $\Gamma$ be a non-trivial finite subgroup of $\mathrm{SU}(2)$. The singular complex surface $\mathbb{C}^{2}/\Gamma$ is termed a \emph{Kleinian singularity}.

\begin{itemize}
    \item Subsection 2.1 revisits the ADE root system correspondence for Kleinian singularities.
    
    \item Subsection 2.2 summarizes structural properties of semi-universal deformations and their simultaneous resolutions.
    
    \item Subsection 2.3 details the hyperK\"ahler ALE-Kleinian singularity relationship: The twistor space of a hyperK\"ahler ALE manifold serves as a simultaneous resolution for deformations of the associated singularity (Proposition \ref{prop_twistor}). We further examine the period map of this twistor space—a fundamental tool for our proofs 
    (Lemma \ref{lem_pull_period} and Corollary \ref{cor_strategy}).
    
    \item Subsection 2.4 establishes a crucial isomorphism in the parameter space $Z\otimes\mathbb{R}^{3}$ (Lemma \ref{lem_rotation}), leveraging Kronheimer's results to enable our proof of Theorem \ref{tm_1} in Section 3.
\end{itemize}

\subsection{Kleinian singularities and root systmes}
\begin{definition}\cite[Chapter III]{BHPV_book}
\label{def_resolution_KLeinian}
A \emph{minimal resolution} of a Kleinian singularity $\mathbb{C}^2 /\Gamma$ is a smooth complex surface $\widetilde{\mathbb{C}^2 / \Gamma}$ together with a proper holomorphic map $\eta : \widetilde{\mathbb{C}^2 / \Gamma} \to \mathbb{C}^2 /\Gamma$ such that:
\begin{itemize}
    \item The restriction, $\eta:  \widetilde{\mathbb{C}^2 / \Gamma} \setminus\eta^{-1}(0) \to\mathbb{C}^2 /\Gamma\setminus{0}$, is an isomorphism,
    \item $\widetilde{\mathbb{C}^2 / \Gamma}$ contains no $(-1)$-curves.
\end{itemize}
The set $\eta^{-1}(0)$ is called the \emph{exceptional locus} of $\widetilde{\mathbb{C}^2 / \Gamma}$.
\end{definition}

Any Kleinian singularity admits a unique minimal resolution, denoted by $\eta: \widetilde{\mathbb{C}^2 / \Gamma} \to \mathbb{C}^2 / \Gamma$. The exceptional locus is a union of rational curves $\eta^{-1}(0) = \Sigma_1 \cup \cdots \cup \Sigma_r$, where each $\Sigma_{i} \simeq \mathbb{P}^{1}$. Topologically, these 2-spheres $\Sigma_{i}$ ($i = 1, \dots, r$) generate the second homology group $H_{2}(\widetilde{\mathbb{C}^2 / \Gamma}; \mathbb{Z})$.

We associate a vertex with each $\Sigma_{i}$ and connect two vertices with an edge if the corresponding curves intersect. This construction yields a diagram $\Delta(\Gamma)$ of $\eta^{-1}(0)$. It is well known that $\Delta(\Gamma)$ is an ADE-type Dynkin diagram. Associating $\Gamma$ with its diagram $\Delta(\Gamma)$ induces a bijection \cite[p.111]{slodowy_LNM}:
\begin{align*}
    \left\{\text{finite subgroups of } \mathrm{SU}(2) \right\}/{\sim} \,\, &\longleftrightarrow \,\, \left\{\text{ADE-type Dynkin diagrams}\right\} \\
    \Gamma \,\, &\mapsto \,\, \Delta(\Gamma).
\end{align*}
Thus, the Kleinian singularity $\mathbb{C}^{2}/\Gamma$ is of type A, D, or E according to whether $\Delta(\Gamma)$ is of type A, D, or E.

Denote by $\mathfrak{h}$ the real Cartan algebra associated with $\Delta(\Gamma)$, and by $\Phi$ the corresponding root system. The vertices of $\Delta(\Gamma)$ correspond to simple roots $\theta_{1}, \dots, \theta_{r} \in \mathfrak{h}^{*}$. Let $\Delta$ denote the set of simple roots, which forms a base of $\Phi$. Under the minimal resolution $\widetilde{\mathbb{C}^{2}/\Gamma}$, the 2-spheres $\Sigma_{i}$ ($i = 1, \dots, r$) correspond bijectively to the simple roots $\theta_{i}$. Since these 2-spheres generate $H_{2}( \widetilde{\mathbb{C}^2 / \Gamma}; \mathbb{Z})$, we obtain an embedding of the root system $\Phi \hookrightarrow H_{2}( \widetilde{\mathbb{C}^2 / \Gamma}; \mathbb{Z})$. Hence there are natural isomorphisms:
\begin{equation*}
    H_{2}\big( \widetilde{\mathbb{C}^2 / \Gamma}; \mathbb{R}\big) \simeq \mathfrak{h}^{*}, \quad 
    H^{2}\big( \widetilde{\mathbb{C}^2 / \Gamma}; \mathbb{C}\big) \simeq \mathfrak{h}^{\mathbb{C}}.
\end{equation*}
The Kleinian singularity $\mathbb{C}^{2}/\Gamma$ is said to be of \emph{type $\Phi$ at the singular point $0$}. This classification is intrinsically local.

\subsection{Deformations of Kleinian singularities}

\begin{definition}\cite[p.7]{HR84}
    A \emph{complex space} is a pair $(X,\mathcal{O}_X)$, where $X$ is a topological space $X$ and $\mathcal{O}_{X}$ is a structure sheaf, which satisfies the following properties: 
    \begin{itemize}
        \item[(i)] Every point $x\in X$ admits a neighborhood $U$ such that $U$ can be represented as $U = D\cap Z$, where $D$ is a domain in $\mathbb{C}^{n}$ for some positive integer $n$ and $Z$ is the zero set of some holomorphic functions $f_{1},\cdots,f_{k}$ on $D$.
        \item[(ii)] And the restriction of the structure sheaf $\mathcal{O}_{X}$ to $U$, denoted by $\mathcal{O}_{X}(U)$, is isomorphic to the restriction of the quotient sheaf $\mathcal{O}_{D}/I$ to $U$, where $I$ is the ideal generated by $f_{1},\cdots,f_{k}$.
    \end{itemize}
Morphisms between complex spaces are called \emph{holomorphic maps}. 
\end{definition}

\begin{definition}\cite[Introduction]{Tjurina_locally-semi-universal}
    A \emph{deformation} of a complex space $X^{0}$ is a flat holomorphic map $\pi:\mathcal{X}\to T$ between two complex spaces $\mathcal{X}$ and $T$, together with an isomorphism $\xi:X^{0}\xrightarrow{\sim} \pi^{-1}(t_{0})$, where $t_{0}$ is a point in $T$. We call $t_{0}$ the \emph{distinguished point}. This data is expressed by the following commutative diagram:
    \begin{equation*}
    \begin{tikzcd}
        X^{0} \arrow[r, "\xi"] \arrow[d] & 
        \mathcal{X} \arrow[d, "\pi"] \\
        \{t_{0}\} \arrow[r, hook, "\subset" above] & 
        T
    \end{tikzcd}
    \end{equation*}
    Here $\pi$ is a flat holomorphic map, and $\xi$ is an embedding.
We say a deformation is a \emph{$\mathbb{C}^{*}$-deformation} if the above diagram admits a $\mathbb{C}^{*}$-equivariant action. This means there exist respectively $\mathbb{C}^{*}$-actions on complex spaces $X^{0},\mathcal{X}$ and $T$ such that $\left\{t_{0}\right\}$ is a fixed point in $T$ and the morphisms $\xi$ and $\pi$ are equivariant. 
\end{definition}

\begin{lemma}\cite[Section 8.7]{slo80}\cite[p.681]{Kron1}
\label{lem_semi_branch} 
There exists a $\mathbb{C}^{*}$-deformation of the Kleinian singularity $\mathbb{C}^{2}/\Gamma$, 
denoted by $f: \mathcal{Y} \to \mathcal{V}$, satisfying:
\begin{itemize}
    \item[(1)] $\mathcal{V}$ is isomorphic to $\mathfrak{h}^{\mathbb{C}}/W$, where 
          $\mathfrak{h}^{\mathbb{C}}$ is the complexification of $\mathfrak{h}$, 
          and the distinguished point $0 \in \mathfrak{h}^{\mathbb{C}}/W$ satisfies 
          $f^{-1}(0) \simeq \mathbb{C}^{2}/\Gamma$.
          
    \item[(2)] The {\rm branch locus} $\mathcal{D} \subset \mathcal{V}$ consists of all points $v$ 
          where the fiber $f^{-1}(v)$ is singular.
          
    \item[(3)] $\mathcal{D}$ equals the image under the quotient map 
          $\mathfrak{h}^{\mathbb{C}} \to \mathfrak{h}^{\mathbb{C}}/W$ of 
          $\bigcup_{\alpha \in \Phi} \ker \alpha$.
\end{itemize}
This deformation $f: \mathcal{Y} \to \mathcal{V}$ is the {\rm semi-universal deformation} 
of $\mathbb{C}^{2}/\Gamma$.
\end{lemma}

\begin{lemma}
\label{lem_affine}
     According to \cite[pp. 23-25]{Kas-Schlessinger_versal-deformation_1972}, each fiber of map $f$ in Lemma \ref{lem_semi_branch} is a complex surface (possibly singular) in the affine space $\mathbb{C}^{3}$.
     Moreover, in \cite[Section 1]{Tyu70}, Tjurina pointed out,
     with the ADE classification of Kleinian singularities stated in our subsection 2.1, the space $\mathcal{Y}$ in the semi-universal deformation $f:\mathcal{Y}\to\mathcal{V}$ is a hypersurface in a vector space $\mathbb{C}^{3}\cross\mathbb{C}^{r}$ given by one of the following equations:

    \begin{align*}
    \label{aln_semiuniversal}
    \mathrm{Type}\,A_{r}:\,\,& xy + z^{r+1} + t_{1}z^{r-1} + t_{2}z^{r-2} + \cdots + t_{r-1}z + t_{r} = 0;\\
    \mathrm{Type}\,D_{r}:\,\,&x^{2} + y^{2}z + z^{r-1} + t_{1}z^{r-2} + \cdots + t_{r-1} + 2t_{r}y = 0;\\
    \mathrm{Type}\,E_{6}:\,\,&x^{2} + y^{3} + z^{4} + t_{1}z^{2} + t_{2}z + t_{3}z + y(t_{4}z^{2} + t_{5}z + t_{6}) =0;\\
    \mathrm{Type}\,E_{7}:\,\,&x^{2} + y^{3} + yz^{3} + y(t_{1}z + t_{2}) + t_{3}z^{4} + \cdots + t_{6}z + t_{7} = 0;\\
    \mathrm{Type}\,E_{8}:\,\,&x^{2} + y^{3} + z^{5} + y(t_{1}z^{3} + \cdots + t_{4}) + t_{5}z^{3} + \cdots + t_{8} = 0,\\
    \end{align*}
    and the projection $f$ is simply the projection to the parameter space $\{(t_{1},\cdots,t_{r})\in\mathbb{C}^{r}\}\simeq\mathcal{V}$.
    This explicit expression will not be used later in the text; we would only use the fact that each fiber of $f$ is in $\mathbb{C}^3$. 
\end{lemma}

We denote by $\Tilde{f}: \Tilde{\mathcal{Y}} \to \mathfrak{h}^{\mathbb{C}}$ the deformation which is the pullback  of the semi-universal deformation $f$ by the natural quotient map $q: \mathfrak{h}^{\mathbb{C}}\to \mathfrak{h}^{\mathbb{C}}/W\simeq\mathcal{V}$, that is,
    \begin{equation}
    \label{gr_liftsemi}
    \begin{tikzcd}
        \widetilde{\mathcal{Y}} \arrow[r] \arrow[d, "\widetilde{f}"'] & 
        \mathcal{Y} \arrow[d, "f"] \\
        \mathfrak{h}^{\mathbb{C}} \arrow[r, "q"] & 
        \mathcal{V}
    \end{tikzcd}
    \end{equation}
Then, we have the explicit description about the singularities in the singular fibers of $\Tilde{f}$:
\begin{lemma}\cite[Lemma 1, p.94]{slo80}
\label{lem_type}
For any $\lambda\in \mathfrak{h}^{\mathbb{C}}$, denote $\Phi_{\lambda} := \{\theta\in\Phi \,|\, \theta(\lambda) = 0\}$, which is a root subsystem of the root system $\Phi$. We have:
    \begin{itemize}
        \item [(1)] The fiber $\Tilde{\mathcal{Y}}_{\lambda}$ has singularities if and only if $\lambda\in\bigcup_{\theta\in\Phi} \ker\ \theta$.
        \item [(2)] Let $\Phi_{\lambda} = \Phi_{\lambda,1}\sqcup\cdots\sqcup\Phi_{\lambda,m}$ be the orthogonal decomposition of $\Phi_{\lambda}$ into irreducible root subsystems. Then the fiber $\Tilde{\mathcal{Y}}_{\lambda}$ has exactly $m$ Kleinian singularities of type $\Phi_{\lambda,1},\cdots,\Phi_{\lambda,m}$ at $m$ points respectively. 
    \end{itemize}
As a consequence, any singular fiber $\Tilde{\mathcal{Y}}_{\lambda}$ is a singular complex surface that has only finite isolated singularities.
\end{lemma}
\qed

\begin{definition}\cite[Section 6. Chapter III]{BHPV_book}
\label{def_mresl}
A \emph{minimal resolution} of a singular complex surface $X$ with only isolated singularities at a finite set $S$ is a smooth complex surface $\widetilde{X}$ together with a proper holomorphic map $\eta : \widetilde{X} \to X$ such that:
\begin{itemize}
    \item The restriction to the smooth part, $\eta: \widetilde{X} \setminus\eta^{-1}(S) \to X\setminus S$, is an isomorphism as complex manifolds.
    \item $\Tilde{X}$ contains no $(-1)$-curves.
\end{itemize}
And we call the pre-image of singularities $\eta^{-1}(S)$ the \emph{exceptional locus} of $\widetilde{X}$.
\end{definition}

\begin{lemma}\cite[Theorem 6.2, Chapter III]{BHPV_book}
\label{lem_uni_mini}
    Every normal surface with finitely many singularities admits a unique minimal resolution. In particular, for each $\lambda\in {\frak h}^{\mathbb{C}}$, the fiber $\mathcal{Y}_{q\left(\lambda\right)}\simeq\Tilde{\mathcal{Y}}_{\lambda}$ is normal since all the singularities are Kleinian by Lemma \ref{lem_type}, and then admits a unique minimal resolution.
\end{lemma}

\begin{lemma}\cite[pp.680-681]{Kron1} \cite[Section 8.7]{slo80}
There exists a commutative diagram of holomorphic maps:
    \begin{equation*}
    \begin{tikzcd}
        \mathfrak{Y} \arrow[r, "\psi"] \arrow[d, "\hat{f}"'] & 
        \mathcal{Y} \arrow[d, "f"] \\
        \mathfrak{h}^{\mathbb{C}} \arrow[r, "q"] & 
        \mathcal{V}
    \end{tikzcd}
    \end{equation*}
    such that: 
    \begin{itemize}
        \item $\hat{f}$ is a holomorphic map between smooth complex manifolds.
        \item The restriction of $\psi$ to each fiber $\mathfrak{Y}_{\lambda}$ of $\hat{f}$ is the minimal resolution of the fiber $\mathcal{Y}_{q(\lambda)}$ of $f$, where $\lambda\in\mathfrak{h}^{\mathbb{C}}$.
    \end{itemize}
We call the map $\hat{f}$ is a \emph{simultaneous resolution} of the semi-universal deformation $f$.
\end{lemma}

\begin{lemma}
\label{lem_num_ration_fiber}
    Any rational curve in the fiber $\mathfrak{Y}_{\lambda}$ of $\hat{f}$, $\lambda\in\mathfrak{h}^{\mathbb{C}}$, must lie in the exception locus of $\mathfrak{Y}_{\lambda}$.
\end{lemma}
\begin{proof}
    By Lemma \ref{lem_affine}, $\mathcal{Y}_{q\left(\lambda\right)}$ is a hypersurface in $\mathbb{C}^{3}$. Since any holomorphic map $\mathbb{P}^{1}\hookrightarrow\mathcal{Y}_{q\left(\lambda\right)}\subset\mathbb{C}^{3}$ must be constant, $\mathcal{Y}_{q\left(\lambda\right)}$ contains no rational curves. Suppose that a rational curve in $\mathfrak{Y}_{\lambda}$ is given by a holomorphic map $\mathbb{P}^{1}\hookrightarrow \mathfrak{Y}_{\lambda}$. Then its composition with the resolution $\mathfrak{Y}_{\lambda}\to\mathcal{Y}_{q\left(\lambda\right)}\subset\mathbb{C}^{3}$ is a constant map. Hence by the definition of minimal resolution, the rational curve must lie in the exceptional locus of $\mathfrak{Y}_{\lambda}$.  
\end{proof}

\begin{corollary}
\label{cor_count}For any $\lambda\in \mathfrak{h}^{\mathbb{C}}$, recall that $\Phi_{\lambda} = \Phi_{1,\lambda}\sqcup\cdots\sqcup\Phi_{m,\lambda}$ represents the type of singularities in Lemma \ref{lem_type}. Then we have the counting of rational curves in $\mathfrak{Y}_{\lambda}$:
    \begin{equation*}
        \#\left\{ \mathrm{rational\,\, curves\,\, in\,\,} \mathfrak{Y}_{\lambda} \right\} = 
        \mathrm{rank}\,\Phi_{\lambda}.
    \end{equation*}
\end{corollary}
\begin{proof}
    By Lemma \ref{lem_num_ration_fiber}, the number of rational curves in $\mathfrak{Y}_{\lambda}$ is the sum of the numbers of rational curves in the exceptional locus of each singularity, hence is given by the sum of the ranks of the types $\Phi_{1,\lambda},\cdots,\Phi_{m,\lambda}$ of these singularities in $\mathcal{Y}_{q\left(\lambda\right)}\simeq \Tilde{\mathcal{Y}}_{\lambda}$. Since the decomposition is orthogonal, we have:
    \begin{equation*}
        \mathrm{rank}\ \Phi_{\lambda} = \mathrm{rank}\ \Phi_{1,\lambda} + \cdots + \mathrm{rank}\ \Phi_{m,\lambda}.
    \end{equation*}
\end{proof}


\subsection{Period map of twistor space}
Recall the notations from Subsections 1.1 and 1.2. For any non-trivial finite subgroup $\Gamma < \mathrm{SU}(2)$, consider $\zeta = (\zeta^1, \zeta^2, \zeta^3) \in (Z \otimes \mathbb{R}^3)^\circ$ as defined in Section 1.1. The hyperK\"ahler manifold $X_\zeta$ carries the following:
\begin{itemize}
    \item Three compatible complex structures $I_\zeta$, $J_\zeta$, $K_\zeta$;
    \item Associated K\"ahler forms $\omega_\zeta^k$ ($k = 1,2,3$);
    \item A twistor space denoted by $\mathrm{Z}(X_\zeta)$ with a canonical holomorphic projection $\pi: \mathrm{Z}(X_\zeta) \to \mathbb{P}^1$.
\end{itemize}

\begin{prop}\cite[pp.691-692]{Kron2}
\label{prop_twistor}
    Related to $\mathrm{Z}(X_{\zeta})$, there is a $\mathbb{C}^{*}$-deformation of $\mathbb{C}^{2}/\Gamma$:
    \begin{equation*}
        \phi: Y \to \mathbb{C}^2.
    \end{equation*}
    where the $\mathbb{C}^{*}$-action on $\mathbb{C}^{2}$ is the natural diagonal action and $\phi^{-1}(0)\simeq\mathbb{C}^{2}/\Gamma$. We define a map $\pi^{s}:{\rm Z}^s\to \mathbb{P}^1$ by:
    \begin{equation*}
    \begin{tikzcd}[column sep=large]
        \left(Y \setminus \phi^{-1}(0)\right) \big/ \mathbb{C}^{*} \arrow[r, equal] \arrow[d, "\phi"'] & 
        \mathrm{Z}^{s} \arrow[d, "\pi^{s}"] \\
        \left(\mathbb{C}^{2}\setminus\{0\}\right) \big/ \mathbb{C}^{*} \arrow[r, equal] & 
        \mathbb{P}^1
    \end{tikzcd}
    \end{equation*}
    Then there exists a commutative diagram below:
    \begin{equation}
    \label{gr_1}
    \begin{tikzcd}
        \mathrm{Z}(X_{\zeta}) \arrow[r, "\chi"] \arrow[d, "\pi"'] & 
        \mathrm{Z}^{s} \arrow[d, "\pi^{s}"] \\
        \mathbb{P}^1 \arrow[r, equal] & 
        \mathbb{P}^1
    \end{tikzcd}
    \end{equation}
    such that the restriction of $\chi$ to each fiber of $\pi$ provides the minimal resolution of corresponding fiber of $\pi^{s}$. 
    We call the projection $\pi:\mathrm{Z}(X_{\zeta})\to \mathbb{P}^{1}$ a \emph{simultaneous resolution} of the map $\pi^{s}$.   
\end{prop}
\qed

As the twistor space of a hyperK\"ahler manifold, $\mathrm{Z}(X_{\zeta})$ carries a twisted vertical 2-form 
$\Omega \in H^{0}\big(\mathrm{Z}(X_{\zeta}),\,\Lambda^{2}T_{F}^{*}\otimes\pi^{*}\mathcal{O}_{\mathbb{P}^{1}}(2)\big)$ 
by \cite{HKLM}. Defined on the affine chart $\mathbb{C} \subset \mathbb{P}^{1} = \mathbb{C} \cup \{\infty\}$, 
it is given by:
\begin{equation*}
  \Omega(u) = \omega^{2}_{\zeta} + \mathrm{i}\omega^{3}_{\zeta} + 2u\,\omega^{1}_{\zeta} - u^{2}\left(\omega^{2}_{\zeta} - \mathrm{i}\omega^{3}_{\zeta}\right),
\end{equation*}
where $u \in \mathbb{C}$ and $T_{F}^{*} := \ker \dd\pi$. Writing $u = u_{1} + \mathrm{i}u_{2}$ with $u_{1}, u_{2} \in \mathbb{R}$, 
we have:
\begin{align*}
  \operatorname{Re} \Omega(u) &= 2u_{1}\omega_{\zeta}^{1} + (1 - u_{1}^{2} + u_{2}^{2})\,\omega_{\zeta}^{2} - 2u_{1}u_{2}\,\omega_{\zeta}^{3}, \\
  \operatorname{Im} \Omega(u) &= 2u_{2}\,\omega_{\zeta}^{1} - 2u_{1}u_{2}\,\omega_{\zeta}^{2} + (1 + u_{1}^{2} - u_{2}^{2})\,\omega_{\zeta}^{3}.
\end{align*}
Note that $\Omega(u)$ is a holomorphic two-form on $(X_{\zeta}, I_{\zeta}^{u})$, where $I_{\zeta}^{u}$ is defined 
in the paragraph containing equation (\ref{eq_iden_S_P}). The cohomology class $[\Omega(u)]$ is called the 
\emph{period} of $(X_{\zeta}, I_{\zeta}^{u})$. 
Taking the cohomology class of $\Omega$ yields, under the identification $H^{2}(X_{\zeta},\mathbb{C}) \simeq \mathfrak{h}^{\mathbb{C}}$, 
an element of $H^{0}\big(\mathbb{P}^{1}, \mathfrak{h}^{\mathbb{C}} \otimes \mathcal{O}_{\mathbb{P}^{1}}(2)\big)$. 
Lifting to $\mathbb{C}^{2}$, we obtain the map:
\begin{equation}
\label{eq_periodmap}
\begin{aligned}
    \widetilde{p_{\zeta}} \colon \mathbb{C}^{2} &\to \mathfrak{h}^{\mathbb{C}}, \\
    (z_{1}, z_{2}) &\mapsto z_{1}^{2}[\omega^{2}_{\zeta} + \mathrm{i}\omega^{3}_{\zeta}] + 2z_{1}z_{2}[\omega^{1}_{\zeta}] - z_{2}^{2}[\omega^{2}_{\zeta} - \mathrm{i}\omega^{3}_{\zeta}].
\end{aligned}
\end{equation}
The composition $p_{\zeta} := q \circ \widetilde{p_{\zeta}} \colon \mathbb{C}^{2} \xrightarrow{\widetilde{p_{\zeta}}} 
\mathfrak{h}^{\mathbb{C}} \xrightarrow{q} \mathcal{V}$ is the \emph{period map} of $\mathrm{Z}(X_{\zeta})$, 
where $\mathcal{V} = \mathfrak{h}^{\mathbb{C}}/W$ is the base of the semi-universal deformation of $\mathbb{C}^{2}/\Gamma$ 
from Subsection 2.2. For simplicity, we also refer to $\widetilde{p_{\zeta}}$ as the period map.

\begin{lemma} \cite[p. 693]{Kron2}
\label{lem_pull_period}
    Recall the semi-universal deformation $f: \mathcal{Y}\to\mathcal{V}$ of $\mathbb{C}^{2}/\Gamma$ in Lemma \ref{lem_semi_branch}. Then the deformation $\phi$ is the pullback of the semi-universal deformation by the period map $p_{\zeta}$, that is,  the following diagram holds:
    \begin{equation*}
    \begin{tikzcd}
        Y \arrow[r] \arrow[d, "\phi"'] & 
        \mathcal{Y} \arrow[d, "f"] \\
        \mathbb{C}^{2} \arrow[r, "p_{\zeta}"] & 
        \mathcal{V}
    \end{tikzcd}
    \end{equation*}
\end{lemma}

The singular fibers of (the pullback of) semi-universal deformation $f$ and the types of their singularities are characterized by Lemma \ref{lem_type}. Combining these with Lemma \ref{lem_pull_period} and Proposition \ref{prop_twistor}, we have:
\begin{corollary}
\label{cor_strategy}
    Recall the holomorphic map $\pi:\mathrm{Z}(X_{\zeta}) \to \mathbb{P}^1$. Given $u\in\mathbb{P}^1$, then:
    \begin{itemize}
        \item[(1)] The complex submanifold $\pi^{-1}(u)$ is isomorphic to the minimal resolution of the fiber $\phi^{-1}(\Tilde{u})$ of $\phi$, where $\Tilde{u}$ is any lifting of $u$ with respect to the natural quotient map $\mathbb{C}^2\setminus0 \to \mathbb{P}^1$.
        \item[(2)] $\pi^{-1}(u)$ contains a rational curve if and only if the period $\widetilde{p_{\zeta}}(\Tilde{u})\in\mathfrak{h}^{\mathbb{C}}$ lies in the kernel of some root in $\Phi$.
        \item[(3)] The number of rational curves in $\pi^{-1}(u)$ is given by $\mathrm{rank}\,\Phi_{\lambda}$, where $\lambda := \widetilde{p_{\zeta}}(\Tilde{u})$.
    \end{itemize}
    
\end{corollary}
\begin{proof}
    The statement (1) is immediately obtained from Proposition \ref{prop_twistor}. By Lemma \ref{lem_pull_period}, the fiber $\phi^{-1}(\Tilde{u})$ is isomorphic to the fiber $f^{-1}(p_{\zeta}(\Tilde{u}))$ hence isomorphic to the fiber $\Tilde{f}^{-1}(\widetilde{p_{\zeta}}(\Tilde{u}))$. Thus, statement (2) is obtained from the Lemma \ref{lem_type}. And statement (3) is a straightforward consequence of Corollary \ref{cor_count}.
\end{proof}

\subsection{An isomorphism between hyperK\"ahler ALEs}

\begin{lemma}\cite[Proposition 4.1]{Kron1}
\label{lem_period_iso}
Given $\zeta = \left(\zeta^{1},\zeta^{2},\zeta^{3}\right)\in\left(Z\otimes\mathbb{R}^{3}\right)^\circ$. There exists a linear isomorphism $$\sigma: Z \xrightarrow{\sim} H^{2}(X_{\zeta}^{\mathrm{Diff}},\mathbb{R})(\simeq \mathfrak{h}),$$ such that $\sigma(\zeta^{k}) = [\omega^{k}_{\zeta}]$, for $k = 1,2,3.$ Consequently, the complexification $\sigma_{\mathbb{C}}: Z \otimes \mathbb{C} \to H^{2}(X_{\zeta}^{\mathrm{Diff}},\mathbb{C})$ is also an isomorphism.
\end{lemma}
\qed

\begin{lemma}
    
\label{lem_rotation}
    Without loss of generality, we assume that $u\neq\infty$ and denote $u = u_{1} + \mathrm{i}u_{2}\in\mathbb{C}\hookrightarrow\mathbb{P}^{1}$. Recall the notations for complex structures: $I_{\zeta}^{u} = aI_{\zeta} + bI_{\zeta} + cI_{\zeta}$. For $\zeta = \left(\zeta^{1},\zeta^{2},\zeta^{3}\right)\in\left(Z\otimes\mathbb{R}^{3}\right)^\circ$, let
    \begin{align*}
     \Tilde{\zeta}^{2} &:= 2u_{1}\zeta^{1} + (1-u_{1}^{2} + u_{2}^{2})\zeta^{2} - 2u_{1}u_{2}\zeta^{3},\\
     \Tilde{\zeta}^{3} &:= 2u_{2}\zeta^{1} - 2u_{1}u_{2}\zeta^{2} + (1 + u_{1}^{2} - u_{2}^{2})\zeta^{3}.
    \end{align*}
Then the complex manifold $\left(X_{\zeta}, I_{\zeta}^{u}\right)$ is isomorphic to the complex manifold $\left(X_{\Tilde{\zeta}}, I_{\Tilde{\zeta}}\right)$, where $\Tilde{\zeta} = \left(\Tilde{\zeta}^{1},\Tilde{\zeta}^{2},\Tilde{\zeta}^{3}\right)$ for any $\Tilde{\zeta}^{1}\in Z$ such that $\Tilde{\zeta}\in\left(Z\otimes\mathbb{R}^{3}\right)^\circ$.

\end{lemma}

\begin{proof}
For the complex manifold $\left(X_{\zeta}, I_{\zeta}^{u}\right)$, there is a holomorphic two form $\Omega(u) := 2u\,\omega_{\zeta}^{1} + \left(1 - u^{2}\right)\omega_{\zeta}^{2} + \mathrm{i}\left(1+u^{2}\right)\omega_{\zeta}^{3}$.  By the isomorphism in Lemma \ref{lem_period_iso}, the period $[\Omega(u)]$ of $\left(X_{\zeta}, I_{\zeta}^{u}\right)$ can also be given by an element in $Z\otimes\mathbb{C}$ as follows:
 \begin{equation*}
     \sigma_{\mathbb{C}}^{-1} \,[\Omega(u)] = \zeta^{2} + \mathrm{i}\zeta^{3} + 2u\zeta^{1} - u^{2}\left(\zeta^{2} - \mathrm{i}\zeta^{3}\right).
 \end{equation*}
By the definitions of $\Tilde{\zeta}^{2}$ and $\Tilde{\zeta}^{3}$, the period of  $(X_{\Tilde{\zeta}}, I_{\Tilde{\zeta}})$ in $Z\otimes\mathbb{C}$ is given by $\Tilde{\zeta}^{2} + \mathrm{i}\Tilde{\zeta}^{2}$ which also equals to $\zeta^{2} + \mathrm{i}\zeta^{3} + 2u\zeta^{1} - u^{2}\left(\zeta^{2} - \mathrm{i}\zeta^{3}\right)$. Hence $\left(X_{\zeta}, I_{\zeta}^{u}\right)$ and $(X_{\Tilde{\zeta}}, I_{\Tilde{\zeta}})$ have the same period. Denote by $u'\in\mathbb{P}^1$ the corresponding point of $(1,0,0)\in\mathbb{S}^2$, and $\Tilde{u}'\in\mathbb{C}^2$ any lifting of $u'$. Then,
\begin{equation*}
    \widetilde{p_{\zeta}}(\Tilde{u}) = \widetilde{p_{\Tilde{\zeta}}}(\Tilde{u}') =: \lambda_{0}.
\end{equation*}
As in the proof of Corollary \ref{cor_strategy}, $\left(X_{\zeta}, I_{\zeta}^{u}\right)$ is the minimal resolution of the fiber $\Tilde{f}^{-1}(\lambda_{0})$ of $\Tilde{f}$. 
Similarly, $(X_{\Tilde{\zeta}}, I_{\Tilde{\zeta}})$ is also the minimal resolution of the fiber $\Tilde{f}^{-1}(\lambda_{0})$ of $\Tilde{f}$.
Since the minimal resolution of a normal surface is unique, see Lemma \ref{lem_uni_mini}, they are isomorphic to each other.

\end{proof}
\begin{remark}
    When $u = \infty$,  $I_{\zeta}^{u} = -I_{\zeta}$ and the complex manifold $\left(X_{\zeta},-I_{\zeta}\right)$ is isomorphic to $\left(X_{\Tilde{\zeta}},I_{\Tilde{\zeta}}\right)$, where $\Tilde{\zeta} := \left(\zeta^1, -\zeta^2, \zeta^3\right)$.
\end{remark}

\section{HyperK\"ahler ALEs as minimal resolution of Kleinian singularities}
This section initiates the proofs of our main theorems. Specifically:
\begin{itemize}
    \item Subsection 3.1 establishes Theorem \ref{tm_1};
    \item Subsection 3.2 addresses case (i) of Theorem \ref{tm_2};
    \item The remaining cases of Theorem \ref{tm_2} are deferred to Section 4.
\end{itemize}

\subsection{Proof of Theorem \ref{tm_1}:}
The core idea involves comparing the periods of $I_{\mathrm{std}}$ and $I_\zeta$.\\

\noindent\textbf{Proof of the "if" part.} 
Suppose $I_{\mathrm{std}} \in \{aI_{\zeta} + bJ_{\zeta} + cK_{\zeta} \mid (a, b, c) \in \mathbb{S}^2\}$ for some $\zeta = (\zeta^1, \zeta^2, \zeta^3) \in (Z \otimes \mathbb{R}^{3})^\circ$. First assume $I_{\mathrm{std}} \simeq I_{\zeta}$, so that $(X_{\zeta}, I_{\zeta}) \simeq \widetilde{\mathbb{C}^{2}/\Gamma}$. If $(\zeta^{2}, \zeta^{3}) \neq 0$, Lemma \ref{lem_period_iso} implies that the cohomology class of the holomorphic 2-form $\omega_{\mathbb{C}} = \omega_{2} + \mathrm{i}\omega_{3}$ on $(X_{\zeta}, I_{\zeta})$ is non-zero:
\[
[\omega_{\mathbb{C}}] = [\omega_{2}] + \mathrm{i}[\omega_{3}] = \sigma(\zeta^{2} + \mathrm{i}\zeta^{3}) \neq 0.
\]
Since $H_{2}(X_{\zeta},\mathbb{Z})$ is generated by the exceptional 2-spheres $\{\Sigma_i\}_{i=1}^r$, there exists some $[\Sigma_{i_{0}}] \in H_{2}(X_{\zeta},\mathbb{Z})$ satisfying $\langle \omega_{\mathbb{C}}, [\Sigma_{i_0}] \rangle \neq 0$. Consequently, $\Sigma_{i_{0}}$ cannot be holomorphic in $(X_{\zeta}, I_{\zeta})$, a contradiction. 
Therefore, when $(X_{\zeta}, I_{\zeta}) \simeq \widetilde{\mathbb{C}^{2}/\Gamma}$, we must have $\zeta = (\zeta^{1}, 0, 0)$. Kronheimer's result confirms that $\zeta = (\zeta^{1}, 0, 0) \in (Z\otimes\mathbb{R}^{3})^\circ$ implies $(X_{\zeta}, I_{\zeta}) \simeq \widetilde{\mathbb{C}^{2}/\Gamma}$. Thus the equivalence holds precisely for such $\zeta$.

Now assume $I_{\mathrm{std}} \simeq I_{\zeta}^{u}$ for some $u \in \mathbb{C} \subset \mathbb{P}^{1}$. Choose $\widetilde{\zeta} = (\widetilde{\zeta}^{1}, \widetilde{\zeta}^{2}, \widetilde{\zeta}^{3})$ as defined in Lemma \ref{lem_rotation}, so that $(X_{\zeta}, I_{\zeta}^{u}) \simeq (X_{\widetilde{\zeta}}, I_{\widetilde{\zeta}})$. Since $I_{\mathrm{std}} \simeq I_{\zeta}^{u}$, we have $(X_{\widetilde{\zeta}}, I_{\widetilde{\zeta}}) \simeq \widetilde{\mathbb{C}^{2}/\Gamma}$. The preceding argument implies $\widetilde{\zeta}^{2} = \widetilde{\zeta}^{3} = 0$, leading to the system:
\begin{equation}
\label{eq_sol_u}
\begin{cases}
     2u_{1}\zeta^{1} + (1 - u_{1}^{2} + u_{2}^{2})\zeta^{2} - 2u_{1}u_{2}\zeta^{3} = 0 \\
     2u_{2}\zeta^{1} - 2u_{1}u_{2}\zeta^{2} + (1 + u_{1}^{2} - u_{2}^{2})\zeta^{3} = 0.
\end{cases}
\end{equation}
These equations force $\mathrm{rank}(\zeta) = 1$. \\

\noindent\textbf{Proof of the "only if" part.}
Given $\zeta \in (Z \otimes \mathbb{R}^{3})^\circ$ with $\mathrm{rank}(\zeta) = 1$, we may write $\zeta = \zeta^1 (x, y, z)$ where $\zeta^{1} \in Z \setminus \bigcup_{\theta \in \Phi} D_{\theta}$ and $(x, y, z) \in \mathbb{R}^{3} \setminus \{0\}$. To find $u$ such that $I_{\mathrm{std}} \simeq I_{\zeta}^{u}$, we solve system \eqref{eq_sol_u}, which holds iff the period of $(X_{\zeta}, I_{\zeta}^{u})$ vanishes. Lifting $u$ to $(z_{1}, z_{2}) \in \mathbb{C}^{2}$, we solve:
\begin{equation*}
    z_{1}^{2}(y\zeta^{1} + \mathrm{i}z\zeta^{1}) + 2z_{1}z_{2}x\zeta^{1} - z_{2}^{2}(y\zeta^{1} - \mathrm{i}z\zeta^{1}) = 0.
\end{equation*}
Since $\zeta^{1} \neq 0$, this reduces to finding non-trivial solutions to:
\begin{equation*}
    z_{1}^{2}(y + \mathrm{i}z) + 2z_{1}z_{2}x - z_{2}^{2}(y - \mathrm{i}z) = 0.
\end{equation*}
Such solutions exist iff there are $(w_{1}, w_{2}, w_{3}) \in \mathbb{C}^{3} \setminus \{0\}$ with $w_{1} w_{3} = w_{2}^{2}$ satisfying:
\begin{equation}
\label{eq_aline}
    w_{1}(y + \mathrm{i}z) + w_{2}(2x) + w_{3}(-y + \mathrm{i}z) = 0.
\end{equation}
Projectively, this requires the line \eqref{eq_aline} in $\mathbb{P}^{2}$ to intersect the conic $w_{1} w_{3} = w_{2}^{2}$. Since this intersection is non-empty, solutions always exist. This completes the proof of Theorem \ref{tm_1}(i).\\

We now complete the proof of Theorem \ref{tm_1}. 

First, note that when $I_{\zeta}^{u} \simeq I_{\mathrm{std}}$, the period of $(X_{\zeta}, I^{u}_{\zeta})$ vanishes. This follows from Lemma \ref{lem_type}, since the exceptional locus has type $\Phi$.

Given $\mathrm{rank}(\zeta) = 1$, the period map $p_{\zeta}$ has image in a one-dimensional subspace of $\mathfrak{h}^{\mathbb{C}}$. By \cite[p.693-p.694]{Kron2}, this image cannot be contained in the branch locus $\mathcal{D}$. Since $\mathcal{D}$ is a finite union of hyperplanes (Lemma \ref{lem_semi_branch}), $\mathrm{Im}\, p_{\zeta}$ intersects $\mathcal{D}$ only at the origin. Thus for any $s \in \mathrm{Im}\, p_{\zeta} \setminus \{0\}$, the fiber $\phi^{-1}(s)$ is smooth by Lemma \ref{lem_pull_period}.

The $\mathbb{C}^{*}$-equivariance of $\phi$ implies that $\phi^{-1}(s)$ and $\phi^{-1}(\gamma s)$ are biholomorphic for $\gamma \in \mathbb{C}^{*}$. Consequently, all complex structures in the hyperK\"ahler family distinct from $I_{\mathrm{std}}$ are equivalent. Moreover, Lemma \ref{lem_affine} ensures that each $\phi^{-1}(s)$ is a smooth affine hypersurface in $\mathbb{C}^{3}$.

Therefore, if $(X_{\zeta}, I^{u}_{\zeta})$ contains rational curves, then necessarily $I_{\zeta}^{u} \simeq I_{\mathrm{std}}$. The precise count of such complex structures (which equals 2) will be established in the next subsection.
\qed

\subsection{Proof of Theorem \ref{tm_2}(i)}
By Theorem \ref{tm_1}, the cardinality of $\{u \in \mathbb{P}^1 \mid (X_{\zeta}, I_{\zeta}^{u}) \text{ contains rational curves}\}$ equals that of $\{u \in \mathbb{P}^1 \mid I_{\zeta}^{u} \simeq I_{\mathrm{std}}\}$. As established in the proof of Theorem \ref{tm_1}, this corresponds to the number of intersection points between the projective line \eqref{eq_aline} and the conic $w_1 w_3 = w_2^2$. 

Elementary algebra shows that tangency occurs iff:
\begin{equation*}
    x^2 + |y + \mathrm{i}z|^2 = 0,
\end{equation*}
which implies $(x, y, z) = (0,0,0)$. This would contradict $\mathrm{rank}(\zeta) = 1$ since $\zeta = \zeta^1(0,0,0) = 0$. Therefore, the line and conic intersect transversely, yielding exactly two intersection points. Consequently:
\begin{equation}
    \#\left\{u \in \mathbb{P}^{1} : \pi^{-1}(u) \text{ contains a rational curve} \right\} = 2.
\end{equation}
Equation \eqref{eq_1_rk1} follows immediately from the $r$ rational curves in $\widetilde{\mathbb{C}^2/\Gamma}$:
\begin{equation}
    \sum_{u \in \mathbb{P}^{1}} \#\left\{\text{rational curves in } \pi^{-1}(u) \right\} = 2r.
\end{equation}
\qed

\medskip

When $\Gamma = \mathbb{Z}_{2}$, we have $\dim Z = 1$, so $\mathrm{rank}(\zeta) = 1$ for all $\zeta \in (Z \otimes \mathbb{R}^{3})^\circ$. Theorems \ref{tm_1} and \ref{tm_2}(i) imply:
\begin{corollary}
\label{cor_Z2}
    In the hyperK\"ahler family of any ALE space asymptotic to $\mathbb{C}^2/\mathbb{Z}_2$, exactly two complex structures are isomorphic to the minimal resolution $\widetilde{\mathbb{C}^2/\mathbb{Z}_2}$.
\end{corollary}
\begin{remark}
    Further properties of the twistor space are discussed in Example \ref{ex_twi_Z2} (Section 5).
\end{remark}

\section{Counting rational curves}
\label{sec_counting}

This section completes the proof of Theorem \ref{tm_2}. Recall the three cases for the parameter $\zeta$: (i) $\operatorname{rank}(\zeta) = 1$; (ii) $\operatorname{rank}(\zeta) = 2$; (3) $\operatorname{rank}(\zeta) = 3$.
Case (i) was established in Section 3. We now prove Cases (ii) and (iii) in Subsections 4.2--4.3, with Corollary \ref{cor_3} addressed in Subsection 4.4.
As a technical preparation for Case (ii), Subsection 4.1 introduces a decomposition of irreducible root systems and presents key combinatorial results of independent interest.

\medskip

\noindent\textbf{Proof strategy overview.} 
By Corollary \ref{cor_strategy} and Lemma \ref{lem_type}, the image of the period map $\widetilde{p_{\zeta}}$ encodes information about:
\begin{enumerate}
    \item Complex structures containing rational curves
    \item The number of such curves
\end{enumerate}
To bound the two cardinalities
\begin{align}
    \label{eq_count_complex_structure}
    &\# \left\{ u \in \mathbb{P}^{1} \mid \pi^{-1}(u) \text{ contains a rational curve} \right\} \\
    \label{eq_count_rational_curve}
    &\sum_{u \in \mathbb{P}^{1}} \# \left\{ \text{rational curves in } \pi^{-1}(u) \right\}
\end{align}
we analyze the geometric position of $\operatorname{Im} \widetilde{p_{\zeta}}$. This position is determined by its intersection with root kernels, which induces a root system decomposition. The decomposition in Subsection 4.1, designed for Case (ii), also provides a framework for Case (iii). We develop distinct approaches for these cases in subsequent subsections.

\subsection{Decomposition of irreducible root systems}
\label{subsec_decomp}

For a nontrivial finite subgroup $\Gamma \leq \mathrm{SU}(2)$, we recall from Subsection 2.1 that $\mathfrak{h}^{\mathbb{C}}$, $\Phi$, and $\Delta$ denote respectively the complexified Cartan subalgebra, the root system, and the set of simple roots. Throughout, we assume $\operatorname{rank} \Phi \geq 2$, i.e., $\Gamma \neq \mathbb{Z}_{2}$. 

Fix a two-dimensional subspace $L \subset \mathfrak{h}^{\mathbb{C}}$ not contained in any root kernel $\ker \theta$ for $\theta \in \Phi$. This induces for each root $\theta \in \Phi \subset \operatorname{Hom}(\mathfrak{h}, \mathbb{R}) \hookrightarrow \operatorname{Hom}(\mathfrak{h}^{\mathbb{C}}, \mathbb{C})$ a one-dimensional kernel $l_\theta := \ker \theta|_L$, yielding the map:
\begin{equation*}
    \begin{aligned}
        \mu_{L} \colon \Phi &\to \mathbb{P}(L) \\
        \theta &\mapsto [l_{\theta}].
    \end{aligned}
\end{equation*}
The image $\operatorname{Im} \mu_L = \{p_1, \dots, p_s\}$ is finite, inducing a disjoint decomposition:
\begin{equation}
\label{dec_rk2}
  \Phi = \bigsqcup_{k=1}^s \Phi_k, \quad
  \Delta = \bigsqcup_{k=1}^s \Delta_k,
\end{equation}
where $\Phi_k := \mu_L^{-1}(p_k)$ and $\Delta_k := \Phi_k \cap \Delta$. Each $\Phi_k$ is a root subsystem satisfying $\Phi_k = \Phi_k \cap \,\operatorname{span}_\mathbb{R} \Phi_k$. Crucially, this decomposition is generally non-orthogonal and $\Delta_k$ may not be a base of $\Phi_k$.
Since $\bigcap_{\theta \in \Phi} \ker \theta = \{0\} \subset \mathfrak{h}$, we have $\bigcap_{\theta \in \Phi} l_\theta = \{0\} \subset L$, implying $s \geq 2$. Moreover:

\begin{lemma}
\label{lem_4_sgeq3}
    $3 \leq s \leq \frac{1}{2} |\Phi|$.
\end{lemma}
\begin{proof}
The upper bound $s \leq \frac{1}{2} |\Phi|$ is immediate. For the lower bound, suppose toward contradiction that $s = 2$. Then $\operatorname{Im} \mu_L = \{[l_{\theta_1}], [l_{\theta_2}]\}$ for distinct roots $\theta_1, \theta_2$, inducing:
\begin{equation*}
    \Phi = \Phi_1 \sqcup \Phi_2, \quad
    \Delta = \Delta_1 \sqcup \Delta_2.
\end{equation*}

Both $\Delta_1$ and $\Delta_2$ are non-empty: If $\Delta_1 = \emptyset$, then $\Delta_2 = \Delta$ implies $\Phi_2 = \Phi$, forcing $\Phi_1 = \emptyset$. Since $\Phi$ is irreducible, there exist $\alpha_1 \in \Delta_1$ and $\alpha_2 \in \Delta_2$ with $(\alpha_1, \alpha_2) \neq 0$ in the ambient Euclidean space. Thus $\alpha_1 \pm \alpha_2 \in \Phi$ for some choice of sign.

Assume $\beta := \alpha_1 + \alpha_2 \in \Phi$. The kernel condition gives:
\begin{align*}
    l_\beta &= \ker(\alpha_1 + \alpha_2)|_L \\
            &= \ker \alpha_1|_L \cap \ker \alpha_2|_L \\
            &= l_{\alpha_1} \cap l_{\alpha_2}.
\end{align*}
Since $\mu_L(\alpha_1) = [l_{\theta_1}] \neq [l_{\theta_2}] = \mu_L(\alpha_2)$, we have $l_{\alpha_1} \neq l_{\alpha_2}$. Thus $\dim(l_{\alpha_1} \cap l_{\alpha_2}) < \dim l_{\alpha_1} = 1$, implying $l_\beta = \{0\}$. But $\beta \in \Phi$ requires $\ker \beta|_L = l_\beta$ to be one-dimensional, contradiction.

Therefore $s \geq 3$ as required.
\end{proof}

Now we define $r := \operatorname{rank} \Phi$, $r_k := \operatorname{rank} \Phi_k$, and $\delta_k := |\Delta_k|$. Note that $\delta_k \leq r_k$ and $r \leq \sum_{k=1}^s r_k$. Moreover:

\begin{lemma}
\label{lem_4_sum_r+1}    
    \begin{equation}
    \label{eq_bestlow}
        r + 1 \leq \sum_{k=1}^s r_k \leq \frac{1}{2} |\Phi|.
    \end{equation}
\end{lemma}

\begin{proof}
The upper bound follows since finer decompositions increase the rank sum. For the lower bound, suppose $r = \sum_{k=1}^s r_k$. As $r = \sum_{k=1}^s \delta_k$ and $\delta_k \leq r_k$, we have $r_k = \delta_k$ for each $k$. Thus $\Delta_k$ forms a base for $\Phi_k$.
Let $\theta_{\max}$ be the maximal root of $\Phi$ relative to $\Delta$, with expression $\theta_{\max} = \sum_{\alpha \in \Delta} k_\alpha \alpha$ ($k_\alpha > 0$). Since each $\Delta_k$ is a proper subset of $\Delta$ (as $s \geq 3$ by Lemma \ref{lem_4_sgeq3}), $\theta_{\max}$ involves roots from multiple $\Delta_k$. Hence $\theta_{\max} \notin \Phi_k$ for any $k$, contradicting $\Phi = \bigsqcup_k \Phi_k$.
\end{proof}

To achieve the stronger lower bound $2r-1$ (improving Lemma \ref{lem_4_sum_r+1}), we introduce a refined class of root system decompositions:

\begin{definition}
A decomposition $\Phi = \bigsqcup_{k=1}^s \Phi_k$ is called \emph{type-1} if:
\begin{enumerate}
    \item $s \geq 1$,
    \item Each $\Phi_k$ is a non-empty proper root subsystem satisfying $\Phi_k = \Phi \cap \operatorname{span}_\mathbb{R} \Phi_k$,
    \item The subsystems are pairwise disjoint.
\end{enumerate}
The space of all type-1 decompositions is denoted $\mathfrak{D}^1(\Phi)$.
\end{definition}

\begin{remark}
While the $L$-induced decomposition \eqref{dec_rk2} is type-1, $\mathfrak{D}^1(\Phi)$ contains decompositions not arising from any subspace $L$. Thus $\mathfrak{D}^1(\Phi)$ provides a coarser classification than the $L$-induced decomposition space.
\end{remark}

For any type-1 decomposition $\Phi = \bigsqcup_{k=1}^s \Phi_k$, define $r_k := \operatorname{rank} \Phi_k$. The \emph{minimal rank sum} $f_1(\Phi)$ is defined as:
\begin{equation*}
    f_1(\Phi) := \inf_{\mathfrak{D}^1(\Phi)} \left\{ \sum_{k=1}^s r_k \right\}.
\end{equation*}
By Lemma \ref{lem_4_sum_r+1}, whose proof depends only on $\Phi$'s irreducibility and properties (1)--(2), we have $f_1(\Phi) \geq r + 1$. A sharper lower bound holds for ADE-type root systems:

\begin{prop}
\label{prop_4_better_lower_rk2}
Let $\Phi(n)$ denote an irreducible ADE root system of rank $n \geq 2$, i.e., $\Phi(n) = A_n$ ($n \geq 2$), $D_n$ ($n \geq 3$), or $E_n$ ($n=6,7,8$). Then:
\begin{equation*}
    2n - 1 \leq f_1(\Phi(n)).
\end{equation*}
\end{prop}

This is proved by induction on $n$. We first establish a key lemma. Consider $\Phi(n+1)$ (an ADE root system of rank $n+1$) embedded in a certain Euclidean space $\mathrm{E}^{n+1}$. Fix an embedding $A_n \hookrightarrow \Phi(n+1)$ (to be specified in the proof). The \emph{set of extra roots} of $\Phi(n+1)$ with respect to $A_n$ is defined as $S := \Phi(n+1) \setminus A_n$.

\begin{lemma}
\label{lem_extra_roots}
The extra roots in $S$ span $\mathrm{E}^{n+1}$, that is, $\dim_{\mathbb{R}} \operatorname{span}_{\mathbb{R}} S = n+1$.
\end{lemma}

\begin{proof}
    The proof involves explicitly calculating for each type of root systems.
    \begin{itemize}
        \item For $\Phi(n+1) = A_{n+1}$.
        Let $\left(\epsilon_{i}\,,\,i = 1,\cdots,n+2\right)$ be the canonical basis of $ \mathbb{R}^{n+2}$,
        and set $\mathrm{E} := \left\{\sum_{i} \epsilon_{i} = 0\right\}\subset \mathbb{R}^{n+1}$. Then by \cite{Lie_Bour_46},    
        the positive roots of $A_{n+1}$ are: $\epsilon_{i}-\epsilon_{j}\,,\,1\leq i<j\leq n+2$. We choose a root base as:
    \begin{align*}
        &\alpha_{1} = \epsilon_{1} - \epsilon_{2}\\
        &\cdots\\
        &\alpha_{i} = \epsilon_{i} - \epsilon_{i+1}\\
        &\cdots\\
        &\alpha_{n} = \epsilon_{n} - \epsilon_{n+1}\\
        &\alpha_{n+1} = \epsilon_{n+1} - \epsilon_{n+2}
    \end{align*}
    Then we can write the positive roots as linear combinations of these simple roots as:
    \begin{align*}
        \epsilon_{i} - \epsilon_{j} &= \sum_{i\leq k<j}\alpha_{k}  \quad ,\,\,1\leq i<j\leq n+2;\\
    \end{align*}
    The first $n$ simple roots $\alpha_{1},\cdots,\alpha_{n}$ determine the embedding of $A_{n}\hookrightarrow A_{n+1}$. The 'extra' positive roots with respect to $A_{n}$ are: $S^{+} := \left\{\epsilon_{i}-\epsilon_{n+2}, 1\leq i\leq n+ 1\right\}$. By elementary computation, we see:
    \begin{equation*}
        \mathrm{span}_{\mathbb{R}}\,S^{+} = \mathrm{E}.
    \end{equation*}
    Thus, $\dim \mathrm{span}_{\mathbb{R}}S = n+1$.

    \item For $\Phi(n+1) = D_{n+1}$:
    We set $\mathrm{E} = \mathbb{R}^{n+1}$ with canonical basis $\left(\epsilon_{i}\,,\,i = 1,\cdots,n+1\right)$.
    Then, by \cite{Lie_Bour_46},
    the roots of $D_{n+1}$ are: $\pm\epsilon_{i}\pm\epsilon_{j}\,,\,1\leq i<j\leq n+1$. We choose a root base as:
    \begin{align*}
        &\alpha_{1} = \epsilon_{1} - \epsilon_{2}\\
        &\cdots\\
        &\alpha_{i} = \epsilon_{i} - \epsilon_{i+1}\\
        &\cdots\\
        &\alpha_{n} = \epsilon_{n} - \epsilon_{n+1}\\
        &\alpha_{n+1} = \epsilon_{n} + \epsilon_{n+1}
    \end{align*}
    Then we draw the Dynkin diagram of $D_{n+1}$ together with an embedding of $A_{n}$ as: 
    
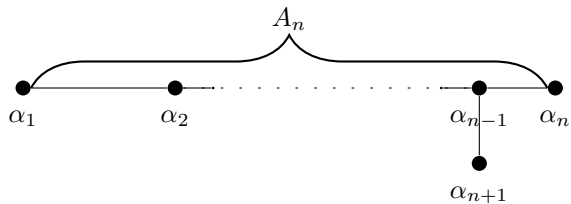
\begin{center}
    \begin{tikzpicture}
    \tikzset{
        node/.style={circle, fill, inner sep=2pt}
        }
    \node[node] (A) at (0,0) {};
        \node at (A.south)[below=2pt] {$\alpha_{1}$};
    \node[node] (B) at (2,0) {};
        \node at (B.south)[below=2pt] {$\alpha_{2}$};
    \node[node] (C) at (6,0) {};
        \node at (C.south)[below=2pt] {$\alpha_{n-1}$};
    \node[node] (D) at (7,0) {};
        \node at (D.south)[below=2pt] {$\alpha_{n}$};
    \node[node] (E) at (6,-1) {};
        \node at (E.south)[below=2pt] {$\alpha_{n+1}$};
    \node[circle, fill, inner sep=0.1pt] (B1) at (2.5,0) {};
    \node[circle, fill, inner sep=0.1pt] (B2) at (5.5,0) {};
    \draw (A) -- (B);
    \draw (B) -- (B1);
    \draw (B2) -- (C);
    \draw[dash pattern=on 1pt off 6pt] (B) -- (C);
    \draw (C) -- (D);
    \draw (C) -- (E);
    \draw[thick, decorate,decoration={brace,amplitude=20pt}]
        (A) -- (D) node [midway, left=10pt, above=18pt] {$A_{n}$};
    
\end{tikzpicture}

    \end{center}

    \noindent We also write the positive roots as linear combinations of these simple roots:
    \begin{align*}
        \epsilon_{i} - \epsilon_{j} &= \sum_{i<k<j}\alpha_{k} & 1\leq i<j\leq n+1;\\
        \epsilon_{i} + \epsilon_{j} &= \sum_{i\leq k<j}\alpha_{k} + 2\left(\sum_{j\leq k<n}\alpha_{k}\right) + \alpha_{n} + \alpha_{n+1} & 1\leq i<j< n+1;\\
        \epsilon_{i} + \epsilon_{n+1} &= \sum_{i\leq k\leq n-1}\alpha_{k} + \alpha_{n+1} & 1\leq i < n+1.
    \end{align*}
    The simple roots $\alpha_{1},\cdots,\alpha_{n}$ determine an embedding of $A_{n}\hookrightarrow  D_{n+1}.$ The 'extra' (positive) roots with respect to $A_{n}$ are: $S^{+} := \left\{\epsilon_{i}+\epsilon_{j}, 1\leq i<j\leq n+ 1\right\}$. By elementary computation, we see:
    \begin{equation*}
        \mathrm{span}_{\mathbb{R}}\,S = \mathrm{E}.
    \end{equation*}
    Thus, $\dim \mathrm{span}_{\mathbb{R}}S = n+1$.

    \item For $\Phi(n+1) = E_{6}, E_{7}, E_{8}$. We only show here the case $\Phi(n+1) = E_{6}$, since the other cases are basically the same.
    The Dynkin diagram of $E_{6}$ together with an embedding of $A_{5}$ is 
        
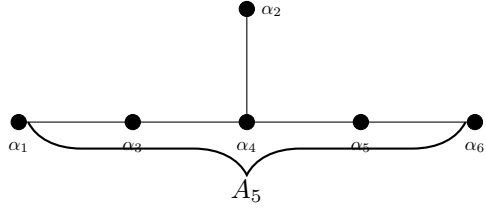
\begin{center}
        \begin{tikzpicture}
        \dynkin[edge length=1.5cm, root radius=0.1cm, labels={\alphalabel{1},\alphalabel{2},\alphalabel{3},\alphalabel{4},\alphalabel{5},\alphalabel{6}}]E6;
        \draw[thick, decorate,decoration={brace,amplitude=20pt,mirror}]
        (root 1) -- (root 6) node [midway, left=10pt, below=18pt] {$A_{5}$};
        \end{tikzpicture}    
        \end{center}

    \noindent where we choose a root base $\left\{\alpha_{1},\cdots,\alpha_{6}\right\}$. The simple roots $\alpha_{1},\alpha_{3},\cdots,\alpha_{6}$ determine an embedding of $A_{5}\hookrightarrow E_{6}$. 
        We denote a root $a\alpha_{1} + b\alpha_{2} + c\alpha_{3} + d\alpha_{4} + e\alpha_{5} + f\alpha
        _{6}$ in the root system $E_{6}$ by a diagram $\begin{pmatrix}
             & & b & & \\
            a&c& d &e&f
        \end{pmatrix}$. Then we cite from \cite{Lie_Bour_46} that the following roots are part of the extra roots of $E_{6}$ with respect to $A_{5}$:
        \begin{align*}
            \begin{pmatrix}
             & & 1 & & \\
            0&1& 2 &1&0
            \end{pmatrix} &
            \begin{pmatrix}
             & & 1 & & \\
            1&1& 2 &1&0
            \end{pmatrix} &
            \begin{pmatrix}
             & & 1 & & \\
            0&1& 2 &1&1
            \end{pmatrix} &
            \begin{pmatrix}
             & & 1 & & \\
            1&2& 2 &1&0
            \end{pmatrix} &
            \begin{pmatrix}
             & & 1 & & \\
            1&1& 2 &1&1
            \end{pmatrix} \\
            \begin{pmatrix}
             & & 1 & & \\
            0&1& 2 &2&1
            \end{pmatrix} &
            \begin{pmatrix}
             & & 1 & & \\
            1&2& 2 &1&1
            \end{pmatrix} &
            \begin{pmatrix}
             & & 1 & & \\
            1&1& 2 &2&1
            \end{pmatrix} &
            \begin{pmatrix}
             & & 1 & & \\
            1&2& 2 &2&1
            \end{pmatrix} &
            \begin{pmatrix}
             & & 1 & & \\
            1&2& 3 &2&1
            \end{pmatrix} \\
            \begin{pmatrix}
             & & 2 & & \\
            1&2& 3 &2&1
            \end{pmatrix} 
        \end{align*}

        \noindent As an observation, the $\mathbb{R}$-span of these 'extra' roots of $E_{6}$ with respect to $A_{5}$ contains all the simple roots $\alpha_{1},\cdots,\alpha_{6}$. For example, the difference of the first two roots above is $\alpha_{1}$, hence $\alpha_{1}$ lies in the span of extra roots. Thus, $\dim \mathrm{span}_{\mathbb{R}}S = 6$.
        
    \end{itemize}
\end{proof}

\noindent\textbf{Proof of Proposition \ref{prop_4_better_lower_rk2}:}
Since $\Phi(n+1)$ is a finite set, the set of all decompositions is also finite. Hence, we can take a decomposition of $\Phi(n+1)$ 
\begin{equation*}
    \Phi(n+1) = \Phi_{1}\sqcup \cdots \sqcup \Phi_{s}\quad ,\quad r_{k} := \mathrm{rank}\,\Phi_{k}.
\end{equation*}
such that $\sum_{k=1}^{s}r_{k} = f_{1}(\Phi(n+1))$. By deleting the extra roots, we obtain the following:

\begin{equation}
\label{dec_lem_extra}
    A_{n} = \Phi_{1}'\sqcup \cdots \sqcup \Phi_{s}' \quad ,\quad \Phi_{k}' := \Phi_{k}\setminus S = \Phi_k\cap A_n\,,\,r_{k}' = \mathrm{rank}\,\Phi_{k}'.
\end{equation}
By definition, each set $\Phi_{k}'$ for $k = 1,\cdots,s$ is either empty or a root subsystem of $A_{n}$. Set $s' := \left\{k \mid \Phi_{k}' \neq \emptyset \right\}$, then we have two possibilities:
\begin{itemize}
    \item $s' = 1$. Then $\Phi_{k_{0}} = A_{n}$ for some $k_{0}$ and the other $\Phi_{k}$, $k\neq k_{0}$, contains only extra roots. By Lemma \ref{lem_extra_roots}:
    \begin{align*}
        f_{1}\left(\Phi(n+1)\right) &= \sum_{k = 1}^{s} r_{k}\\
         & = r_{k_{0}} + \sum_{k\neq k_{0}} r_{k}\\
         & \geq n + n+1 = 2n+1.
    \end{align*}
    \item $s'>1$. Then the disjoint union (\ref{dec_lem_extra}) is a decomposition of $A_{n}$ of type-1 (ignoring the empty sets). By Lemma \ref{lem_extra_roots}, these extra roots cannot be put into a single root subsystem, hence there exist at least two root subsystems $\Phi_{k_{i}}, i = 1,2$ such that $r_{k_{i}} \geq r_{k_{i}}' + 1$. Therefore, 
    \begin{align*}
        f_{1}\left(\Phi(n+1)\right) & = \sum_{k=1}^{s} r_{k}\\
         &\geq \sum_{k = 1}^{s} r_{k}' + 2 \\
         & \geq f_{1}\left(A_{n}\right) + 2.
    \end{align*}
\end{itemize}
Now we can prove the lower bound $2n-1$ for $A_{n}$ in the proposition by induction on $n$:

Lemma \ref{lem_4_sum_r+1} implies that $f_{1}\left(A_{2}\right)\geq 3 = 2\cdot2-1$. And the assumption $f_{1}\left(A_{n}\right) \geq 2n-1$ implies that either $f_{1}\left(A_{n+1}\right) \geq 2n+1$ or $f_{1}\left(A_{n+1}\right) \geq f_{1}\left(A_{n}\right) + 2 \geq 2n-1 + 2 = 2n+1$. Thus, by induction we have $f_{1}\left(A_{n}\right)\geq 2n-1$. And then it is a repeated argument to show $f_{1}\left(\Phi_{n}\right)\geq 2n-1$.
\qed

\subsection{Proof of Theorem \ref{tm_2}(ii)} 
Assume $\operatorname{rank}(\zeta) = 2$. Without loss of generality, $\zeta^{1} \in \operatorname{span}\{\zeta^{2}, \zeta^{3}\}$. The period map:
\begin{equation*}
    \widetilde{p_{\zeta}} \colon \mathbb{C}^{2} \to \mathfrak{h}^{\mathbb{C}}, \quad
    (z_{1}, z_{2}) \mapsto 
    z_{1}^{2}(\zeta^{2} + \mathrm{i}\zeta^{3}) + 2z_{1}z_{2}\zeta^{1} - z_{2}^{2}(\zeta^{2} - \mathrm{i}\zeta^{3})
\end{equation*}
factors through $\mathbb{P}^1$ via:
\begin{align*}
    p_{\zeta}' \colon \mathbb{P}^{1} &\to \mathbb{P}(\operatorname{Im} \widetilde{p_{\zeta}}) \\
    [z_{1} : z_{2}] &\mapsto [w_{1} : w_{2}] := [z_{1}^{2} + 2cz_{1}z_{2} : -z_{2}^{2} + 2\bar{c}z_{1}z_{2}]
\end{align*}
where $\zeta^{1} = c(\zeta^{2} + \mathrm{i}\zeta^{3}) + \bar{c}(\zeta^{2} - \mathrm{i}\zeta^{3})$ for some $c \in \mathbb{C}$. This is a degree-2 branched covering with exactly two branch points.

Since $\operatorname{rank}(\zeta) = 2$, $\operatorname{Im} \widetilde{p_{\zeta}}$ is a 2-dimensional subspace satisfying the $L$-condition from Subsection 4.1 (by \cite[p.693-694]{Kron2}). Setting $L := \operatorname{Im} \widetilde{p_{\zeta}}$, we obtain the decomposition:
\begin{equation*}
    \Phi = \bigsqcup_{k=1}^s \Phi_k, \quad
    \Delta = \bigsqcup_{k=1}^s \Delta_k
\end{equation*}
with $s \geq 3$ and $2r-1 \leq \sum_{k=1}^s r_k \leq \frac{1}{2} |\Phi|$. Let $\{l_k\}_{k=1}^s$ be the corresponding lines in $L$.

By Corollary \ref{cor_strategy}:
\begin{equation}
    \# \left\{ u \in \mathbb{P}^{1} : \pi^{-1}(u) \text{ contains a rational curve} \right\} 
    = \# (p_{\zeta}')^{-1} \big( \{ [l_1], \dots, [l_s] \} \big).
\end{equation}
The upper bound $|\Phi|$ is known \cite[Lemma 2.7]{Kron2}. For the lower bound, note that each $[l_k]$ has two preimages $\{u_k^1, u_k^2\} \subset \mathbb{P}^1$, coinciding iff $[l_k] \in \{[1:0], [0:1]\}$. Since $s \geq 3$:
\begin{equation}
    \# (p_{\zeta}')^{-1} \big( \{ [l_1], \dots, [l_s] \} \big) \geq 4.
\end{equation}
Each $(X_{\zeta}, u_k^i)$ ($i=1,2$) is the minimal resolution of $\widetilde{f}^{-1}(v)$ for $v \in l_k \setminus \{0\}$ (see Diagram \ref{gr_liftsemi}). By Corollary \ref{cor_count}, $\widetilde{f}^{-1}(v)$ contains $r_k$ rational curves. Proposition \ref{prop_4_better_lower_rk2} yields:
\begin{align*}
    2r \leq \sum_{k=1}^s r_k \cdot \# \{u_k^1, u_k^2\} 
    &= \sum_{u \in \mathbb{P}^1} \# \{\text{rational curves in } \pi^{-1}(u)\} 
    \leq |\Phi|.
\end{align*}
This completes the proof for Case (ii).

\subsection{Proof of Theorem \ref{tm_2}(iii)} 
Assume $\operatorname{rank}(\zeta) = 3$. The proof strategy begins with analyzing the period map's image. Since $\dim V_\zeta = 3$ where $V_\zeta := \operatorname{span}_{\mathbb{C}} \{ \zeta^1, \zeta^2, \zeta^3 \} \subset \mathfrak{h}^{\mathbb{C}}$, the image $\operatorname{Im} \widetilde{p_\zeta}$ is a hypersurface in $V_\zeta$. 

Choosing the basis $\{ \zeta^2 + i\zeta^3,\, 2\zeta^1,\, -\zeta^2 + i\zeta^3 \}$ for $V_\zeta$, the period map induces a bijection:
\begin{align*}
    p_{\zeta}' \colon \mathbb{P}^1 &\xrightarrow{\sim} C \subset \mathbb{P}(V_\zeta) \simeq \mathbb{P}^2 \\
    [z_1:z_2] &\mapsto [z_1^2 : z_1z_2 : z_2^2]
\end{align*}
Here $C = \mathbb{P}(\operatorname{Im} \widetilde{p_\zeta})$ is the smooth conic $\{xy = z^2\} \subset \mathbb{P}^2$. 

As $\operatorname{Im} \widetilde{p_\zeta} \not\subset \ker \theta$ for any $\theta \in \Phi$, define for each root: $H_\theta := \ker \theta|_{V_\zeta}$ (a 2-dimensional subspace); $\mathbb{P}H_\theta \simeq \mathbb{P}^1$ (projective line in $\mathbb{P}(V_\zeta)$).
Each $\mathbb{P}H_\theta$ intersects $C$ either transversely (at two points) or tangentially (at one point). Hence, for any intersection point $\gamma\in C$ we define
\begin{definition}
    An intersection point $\gamma$ is called a \emph{tangential point} if all $\mathbb{P}H_\theta$ through $\gamma$ are tangent to $C$;
    $\gamma $  is called a \emph{transversal point} if some $\mathbb{P}H_\theta$ intersects $C$ transversely at $\gamma$.
\end{definition}

The intersection locus decomposes as:
\[
\bigcup_{\theta \in \Phi} (\mathbb{P}H_\theta \cap C) = \{\gamma_1, \dots, \gamma_s, \gamma_{s+1}, \dots, \gamma_{s+t}\}
\]
where $\{\gamma_k\}_{k=1}^s$ are transversal points and $\{\gamma_{s+j}\}_{j=1}^t$ are tangential points. 
By the bijectivity of $p_\zeta'$, each $\gamma_l$ for $l = 1,\cdots,s+t$ corresponds to a unique complex structure $u_l \in \mathbb{P}^1$. Thus:
\begin{equation}
\label{eq_root_count_complex_structure_rk3}
s + t = \# \left\{ u \in \mathbb{P}^1 : \pi^{-1}(u) \text{ contains a rational curve} \right\}
\end{equation}

This induces the decomposition:
\begin{equation}
\label{dec_rk3}
\Phi = \bigcup_{l=1}^{s+t} \Phi_l, \quad
\Delta = \bigcup_{l=1}^{s+t} \Delta_l
\end{equation}
where for $l = 1,\dots,s+t$:
\begin{align*}
\Phi_l := \{\theta \in \Phi : \gamma_l \in \mathbb{P}H_\theta \cap C\}, \quad \Delta_l := \Phi_l \cap \Delta
\end{align*}
Each $\Phi_l$ for $l = 1,\cdots,s+t$ is a root subsystem satisfying 
\begin{equation}
\label{cond_root_span}
\Phi \cap \operatorname{span}_{\mathbb{R}} \Phi_l = \Phi_l.
\end{equation}

\begin{lemma}
\label{lem_4_disj}
The subsystems $\{\Phi_{s+j}\}_{j=1}^t$ are pairwise disjoint, and $\bigcup_{k=1}^s \Phi_k$ is disjoint from $\bigcup_{j=1}^t \Phi_{s+j}$.
\end{lemma}

\begin{proof}
For any tangential point $\gamma_{s+j}$, every line $\mathbb{P}H_\theta$ through it is tangent to $C$. As $C$ is a conic, such lines intersect $C$ only at $\gamma_{s+j}$. Thus:
\begin{itemize}
    \item $\Phi_{s+j} \cap \Phi_{s+j'} = \emptyset$ for $j \neq j'$ (distinct tangential points)
    \item $\Phi_{s+j} \cap \Phi_k = \emptyset$ (no line through $\gamma_{s+j}$ passes through $\gamma_k$)
\end{itemize}
\end{proof}

Denote $r_k := \operatorname{rank} \Phi_k$, $\delta_k := |\Delta_k|$, $r_{s+j} := \operatorname{rank} \Phi_{s+j}$, and $\delta_{s+j} := |\Delta_{s+j}|$. Equation \eqref{eq_2_rk3} then becomes:
\begin{equation}
\label{eq_root_count_rational_curve_rk3}
    \sum_{k=1}^{s} r_k + \sum_{j=1}^{t} r_{s+j} = \sum_{u \in \mathbb{P}^{1}} \# \left\{ \text{rational curves in } \pi^{-1}(u) \right\}.
\end{equation}

The upper bound follows by refining the decomposition \eqref{dec_rk3} into disjoint subsystems: 
First, 
\begin{equation*}
    \Phi_1 \cup \cdots \cup \Phi_s = \left(\bigsqcup_{k=1}^s \Phi_{kk} \right)  \sqcup \left(\bigsqcup_{1 \leq k_1 < k_2 \leq s} \Phi_{k_1k_2}\right)
\end{equation*}
where
\begin{align*}
    \Phi_{kk} 
    &:= \left\{ \theta \in \Phi : \mathbb{P}H_\theta \text{ tangent to } C \text{ at } \gamma_k \right\} \subset \Phi_k
    & (k &= 1,\dots,s) \\
    \Phi_{k_1k_2} &:= \left\{ \theta \in \Phi : \mathbb{P}H_\theta \ni \gamma_{k_1}, \gamma_{k_2} \right\} = \Phi_{k_1} \cap \Phi_{k_2} & (1 &\leq k_1 \neq k_2 \leq s).
\end{align*}
Since $\Phi_k = \Phi_{kk}\cup \bigcup_{k'\neq k}\Phi_{kk'}$, we see $r_k \leq r_{kk} + \sum_{k' \neq k} r_{kk'}$, implying $\sum_{k=1}^s r_k \leq  \sum_{k_1 , k_2} r_{k_1k_2}$. Consequently:
\begin{align*}
    \sum_{k=1}^s r_k + \sum_{j=1}^t r_{s+j} 
    &\leq  \sum_k r_{kk} +  \sum_{k_1 \neq k_2} r_{k_1k_2} + \sum_{j=1}^t r_{s+j} \\
    &\leq 2 \left( \sum_{k = 1}^s r_{kk} + \sum_{1\leq k_1 < k_2\leq s} r_{k_1k_2} + \sum_{j=1}^t r_{s+j} \right) \\
    &\leq |\Phi|
\end{align*}
where the last inequality follows from Lemma \ref{lem_4_sum_r+1} since the union $\bigsqcup_{k = 1}^s\Phi_{kk}\sqcup\bigsqcup_{1\leq k_1 < k_2 \leq s} \Phi_{k_1k_2} \sqcup \bigsqcup_{j=1}^t \Phi_{s+j}$ forms a type-1 decomposition of $\Phi$ (implied by Lemma \ref{lem_4_disj}). 

To establish the lower bound $s + t \geq 3$ for \eqref{eq_root_count_complex_structure_rk3}, we consider three cases:
\begin{enumerate}
    \item $s = 0$ (all tangential points)
    \item $t = 0$ (all transversal points)
    \item $s > 0$ and $t > 0$ (mixed types)
\end{enumerate}
For lower bound of \eqref{eq_root_count_rational_curve_rk3}, we introduce in Subsection 4.3 a refined decomposition of $\Phi$ analogous to that in Subsection 4.1 to prove $2r-1 \leq \sum_{l = 1}^{s+t} r_l$ .

\subsubsection{Case $s = 0$: Impossibility} 
We prove this case cannot occur.

\begin{lemma}
\label{lem_4_rk3_case1}
    In decomposition \eqref{dec_rk3}, $s \neq 0$.
\end{lemma}

\begin{proof}
    Suppose $s = 0$. By Lemma \ref{lem_4_disj}, the decomposition simplifies to:
    \begin{equation*}
        \Phi = \bigsqcup_{j=1}^t \Phi_j, \quad
        \Delta = \bigsqcup_{j=1}^t \Delta_j
    \end{equation*}
    where we re-index $\Phi_{0+j}$ as $\Phi_j$. Lemma \ref{lem_4_sgeq3} implies $t \geq 3$. By irreducibility, we may assume $\Phi_1 \neq \emptyset$, $\Phi_2 \neq \emptyset$, and $(\Phi_1, \Phi_2) \neq 0$. Thus there exist $\alpha \in \Phi_1$, $\beta \in \Phi_2$ with $\alpha \pm \beta \in \Phi$. Assume $\gamma := \alpha + \beta \in \Phi_3$ without loss of generality.
    
    The projective intersections satisfy:
    \begin{equation*}
        \mathbb{P}H_\alpha \cap \mathbb{P}H_\beta 
        = \mathbb{P}(H_\alpha \cap H_\beta)
        = \mathbb{P}H_\alpha \cap \mathbb{P}H_\gamma 
        = \mathbb{P}H_\beta \cap \mathbb{P}H_\gamma
    \end{equation*}
    Thus $\mathbb{P}H_\alpha$, $\mathbb{P}H_\beta$, $\mathbb{P}H_\gamma$ are concurrent at some $p \in \mathbb{P}(V_\zeta)$. By definition, each is tangent to $C$ at distinct points $\gamma_1, \gamma_2, \gamma_3$. But a conic cannot have three distinct tangent lines concurrent, a contradiction.
\end{proof}

\subsubsection{Case $t = 0$: Transversal points only}
\begin{lemma}
\label{lem_4_rk3_case2}
    For decomposition \eqref{dec_rk3} with $t = 0$, we have $s \geq 3$.
\end{lemma}

\begin{proof}
    Assume $t = 0$. Since $s$ counts transversal points, minimality implies $s \geq 2$. Suppose toward contradiction $s = 2$, so all intersection points are $\{\gamma_1, \gamma_2\}$. The decomposition becomes:
    \begin{align*}
        \Phi &= \Phi_1 \cup \Phi_2, \\
        \Delta &= \Delta_1 \cup \Delta_2.
    \end{align*}
    Refine this into disjoint subsystems:
    \begin{align*}
        \Phi &= \Phi'_1 \sqcup \Phi'_2 \sqcup \Phi'_3, \\
        \Delta &= \Delta'_1 \sqcup \Delta'_2 \sqcup \Delta'_3,
    \end{align*}
    where
    \begin{align*}
        \Phi'_1 &:= \{\theta : \mathbb{P}H_\theta \text{ tangent to } C \text{ at } \gamma_1\}\\
        \Phi'_2 &:= \{\theta : \mathbb{P}H_\theta \text{ tangent to } C \text{ at } \gamma_2\}\\
        \Phi'_3 &:= \{\theta : \mathbb{P}H_\theta \text{ contains } \gamma_1 \text{ and } \gamma_2\}
    \end{align*}
    These are root subsystems with $\Phi'_3 \neq \emptyset$ and $\Phi'_1, \Phi'_2 \neq \emptyset$ (since $\bigcap_{\theta} H_\theta = \{0\}$). By irreducibility, they are not pairwise orthogonal. Without loss of generality, take $\alpha \in \Phi'_1$, $\beta \in \Phi'_2$ with $\gamma := \alpha + \beta \in \Phi'_3$.

    Projectively, there holds
    \[
        \mathbb{P}H_\alpha \cap \mathbb{P}H_\gamma = \mathbb{P}H_\beta \cap \mathbb{P}H_\gamma = \mathbb{P}H_\alpha \cap \mathbb{P}H_\beta
    \]
    Thus $\mathbb{P}H_\alpha$, $\mathbb{P}H_\beta$, $\mathbb{P}H_\gamma$ are concurrent at some $p \in \mathbb{P}(V_\zeta)$. But we see
    \begin{align*}
        &\mathbb{P}H_\alpha \text{ tangent to } C \text{ at } \gamma_1\\
        &\mathbb{P}H_\beta \text{ tangent to } C \text{ at } \gamma_2\\
        &\mathbb{P}H_\gamma \text{ transverse to } C \text{ at both } \gamma_1 \text{ and } \gamma_2
    \end{align*}
    This configuration is geometrically impossible in $\mathbb{P}^2$. Contradiction.
\end{proof}

Thus for $t = 0$:
\begin{equation}
    3 \leq s = \# \left\{ u \in \mathbb{P}^1 : \pi^{-1}(u) \text{ contains a rational curve} \right\}
\end{equation}


\subsubsection{Case $s > 0$ and $t > 0$: Mixed intersection types}
\label{subsec:mixed_case}

We have $s \geq 2$ and $t \geq 1$, so:
\begin{equation*}
    s + t \geq 3 = \# \{ u_k, u_{s+j} : k=1,\dots,s;\, j=1,\dots,t \} 
    = \# \left\{ u \in \mathbb{P}^1 : \pi^{-1}(u) \text{ contains a rational curve} \right\}.
\end{equation*}

Moreover:

\begin{lemma}
\label{lem_4_rk3_case3}
    For decomposition \eqref{dec_rk3} with $s=2$, we have $t \geq 2$.
\end{lemma}

\begin{proof}
    Assume $s = 2$. By Lemma \ref{lem_4_rk3_case2}, $t \geq 1$. 

    \textbf{Claim:} $\Phi_1 \cup \Phi_2$ is a root subsystem. 
    Suppose not. Then there exist $\alpha, \beta \in \Phi_1 \cup \Phi_2$ and an integer $m$ such that $\gamma := \alpha + m\beta \in \Phi$ but $\gamma \notin \Phi_1 \cup \Phi_2$. Thus $\gamma \in \Phi_{2+j_0}$ for some $j_0$. Since $\Phi_1 \cap \Phi_2$ is a root subsystem, $\alpha$ and $\beta$ cannot both lie in $\Phi_1 \cap \Phi_2$. By symmetry, consider two cases:
    \begin{enumerate}[label=(P\arabic*)]
        \item $\alpha \in \Phi_1 \setminus \Phi_2$, $\beta \in \Phi_2 \setminus \Phi_1$. Then $\mathbb{P}H_\alpha$, $\mathbb{P}H_\beta$, and $\mathbb{P}H_\gamma$ are concurrent. But
        \begin{align*}
            &\mathbb{P}H_\alpha \text{ tangent to } C \text{ at } \gamma_1\\
            &\mathbb{P}H_\beta  \text{ tangent to } C \text{ at } \gamma_2 \\
            &\mathbb{P}H_\gamma \text{ tangent to C at } \gamma_{2+j_0} (\text{ since} \gamma \in \Phi_{2+j_0})
        \end{align*}
            Three distinct tangent lines cannot be concurrent in $\mathbb{P}^2$—contradiction.
        
        \item $\alpha \in \Phi_1 \cap \Phi_2$, $\beta \in \Phi_2 \setminus \Phi_1$. Then
        \begin{align*}
            &
            \mathbb{P}H_\alpha \text{ contains both } \gamma_1 \text{ and } \gamma_2 (\text{ since } \alpha \in \Phi_1 \cap \Phi_2)\\
            &
            \mathbb{P}H_\beta \text{ tangent at }\gamma_2\\
            &
            \mathbb{P}H_\gamma  \text{ must contain } \gamma_2 (\text{ since } \gamma_2 \text{ is transversal })
        \end{align*}
            But $\gamma \in \Phi_{2+j_0}$ implies $\mathbb{P}H_\gamma$ is tangent at $\gamma_{2+j_0} \neq \gamma_2$, contradiction.
    \end{enumerate}
    Thus $\Phi_1 \cup \Phi_2$ is a root subsystem satisfying:
    \begin{equation}
    \label{cond_root_union_span}
        \Phi \cap \operatorname{span}_{\mathbb{R}} (\Phi_1 \cup \Phi_2) = \Phi_1 \cup \Phi_2.
    \end{equation}
    
    \textbf{Proof of $t \geq 2$:} 
    Suppose $t = 1$. Then $\Phi = (\Phi_1 \cup \Phi_2) \sqcup \Phi_3$. By irreducibility, $(\Phi_1 \cup \Phi_2, \Phi_3) \neq 0$. Take $\alpha \in \Phi_1 \cup \Phi_2$, $\theta \in \Phi_3$ with $(\alpha, \theta) \neq 0$. Assume $\eta := \alpha + \theta \in \Phi$. Conditions \eqref{cond_root_span} and \eqref{cond_root_union_span} imply $\eta \notin \Phi_1 \cup \Phi_2 \cup \Phi_3$, contradiction.
\end{proof}

Thus for $\operatorname{rank}(\zeta) = 3$:
\begin{equation}
    3 \leq \# \left\{ u \in \mathbb{P}^1 : \pi^{-1}(u) \text{ contains a rational curve} \right\}.
\end{equation}
We now establish the lower bound for \eqref{eq_2_rk3} or \eqref{eq_root_count_rational_curve_rk3}:
\begin{equation}
    2r-1 \leq \sum_{k=1}^{s} r_k + \sum_{j=1}^{t} r_{s+j} \leq \sum_{u\in\mathbb{P}^1} \# \left\{ \text{rational curves in } \pi^{-1}(u) \right\}.
\end{equation}

To achieve $2r-1$, we introduce a refined decomposition class:

\begin{definition}
A \emph{type-2 decomposition} of $\Phi$ is an expression
\[
\Phi = \left( \bigcup_{k=1}^s \Phi_k \right) \sqcup \left( \bigsqcup_{j=1}^t \Phi_{s+j} \right)
\]
satisfying:
\begin{enumerate}
    \item $s, t \geq 0$ with $s + t \geq 2$
    \item Each $\Phi_k$ is a non-empty root subsystem such that $\Phi_k = \Phi \cap \operatorname{span}_\mathbb{R} \Phi_k$.
    \item The $\Phi_{s+j}$'s are pairwise disjoint, and $(\bigcup_{k=1}^s \Phi_k) \cap (\bigcup_{j=1}^t \Phi_{s+j}) = \emptyset$
\end{enumerate}
The space of all type-2 decompositions is denoted by $\mathfrak{D}^2(\Phi)$.
\end{definition}

For any type-2 decomposition of $ \mathfrak{D}^2(\Phi)$, define $r_l := \operatorname{rank} \Phi_l$ and the \emph{minimal rank sum}:
\begin{equation*}
    f_2(\Phi) := \inf_{\mathfrak{D}^2(\Phi)} \left\{ \sum_{l=1}^{s+t} r_l \right\}.
\end{equation*}
Similar to Lemma \ref{lem_4_sum_r+1}, we have $f_2(\Phi) \geq r + 1$.

\begin{remark}
The conic-induced decomposition \eqref{dec_rk3} is type-2, so $f_2(\Phi)$ bounds \eqref{eq_root_count_rational_curve_rk3} or \eqref{eq_2_rk3}. When $s=0$, type-1 decompositions embed into $\mathfrak{D}^2(\Phi)$. However, $\mathfrak{D}^2(\Phi)$ is too coarse for obtaining the optimal bounds.
\end{remark}

\begin{prop}
\label{prop_4_better_lower_rk33}
For irreducible ADE root systems $\Phi(n)$ of rank $n \geq 2$:
\begin{equation*}
    2n - 1 \leq f_2(\Phi(n)).
\end{equation*}
\end{prop}

\begin{proof}
The proof is basically the same as Proposition \ref{prop_4_better_lower_rk2}. Since the sum $\sum_{l=1}^{s+t} r_l$ is a positive integer, the infimum is attained. Fix a minimal type-2 decomposition:
\[
\Phi(n+1) = \left( \bigcup_{k=1}^s \Phi_k \right) \sqcup \left( \bigsqcup_{j=1}^t \Phi_{s+j} \right).
\]
Let $S$ be the extra roots in $\Phi(n+1)$ relative to $A_n$ (as in Lemma \ref{lem_extra_roots}). By deleting $S$ from $\Phi(n+1)$, we obtain a decomposition for $A_n$:
\begin{equation}
\label{dec_extra_rk3}
A_n = \left( \bigcup_{k=1}^s \Phi'_k \right) \sqcup \left( \bigsqcup_{j=1}^t \Phi'_{s+j} \right), \quad
\Phi'_l := \Phi_l \setminus S = \Phi_l\cap A_n.
\end{equation}
Each $\Phi'_l$ is either empty or a root subsystem. Two cases arise:

\noindent\textbf{Case 1:} Some $\Phi'_{l_0} = A_n$. Then:
\[
f_2(\Phi(n+1)) \geq \operatorname{rank} A_n + \sum_{l \neq l_0} \operatorname{rank} \Phi_l \geq n + (n+1) = 2n+1.
\]

\noindent\textbf{Case 2:} All $\Phi'_l$ are proper subsystems. Then \eqref{dec_extra_rk3} induces a type-2 decomposition of $A_n$ (ignoring empty set). Since $S$ spans $\mathbb{E}^{n+1}$:
\[
f_2(\Phi(n+1)) \geq f_2(A_n) + 2.
\]

\noindent\textbf{Induction:} For $A_3$ there holds $ f_2(A_3) \geq 3 = 2\cdot2-1$. Assume $f_2(A_n) \geq 2n-1$. Then:
\[
f_2(A_{n+1}) \geq \min(2n+1, (2n-1) + 2) = 2n+1.
\]
Thus $f_2(A_n) \geq 2n-1$. The analogous arguments extend to $D_n$ and $E_n$.
\end{proof}

Therefore:
\begin{equation}
    2r-1 \leq f_2(\Phi) \leq \sum_{k=1}^{s} r_k + \sum_{j=1}^{t} r_{s+j} \leq \sum_{u\in\mathbb{P}^1} \# \left\{ \text{rational curves in } \pi^{-1}(u) \right\}.
\end{equation}
This completes the proof of Theorem \ref{tm_2}(iii).

\newpage 
\section{The semi-continuity of $Q_1$}

Recall the map $Q_1$ defined in the introduction: 
\begin{align*}
    Q_1 : \left(Z\otimes\mathbb{R}^3\right)^\circ &\to \mathbb{Z}\\
    \zeta &\mapsto \#\left\{u \in \mathbb{P}^{1} \,|\, \pi_\zeta^{-1}(u) \mathrm{\,\,contains\,\, a\,\, rational\,\, curve} \right\}
\end{align*}
In this section we prove Theorem \ref{tm_3}, which states that
    \begin{equation*}
        \liminf_{\zeta \to \zeta_0} Q_1(\zeta) \geq Q_1(\zeta_0)\quad ,\,\,\forall \zeta_0\in\left(Z\otimes\mathbb{R}^3\right)^\circ.
    \end{equation*}
\begin{proof} We proceed by three cases with respect to $\mathrm{rank}(\zeta_0)$. The proof arise from studying how the plane $L_\zeta$ and the conic $C_\zeta$ behave under perturbation of $\zeta$:
    \begin{itemize}
        \item[(1)] $\text{rank}(\zeta_0) = 1$.
        Theorem \ref{tm_2} implies that for all $\zeta\in\left(Z\otimes\mathbb{R}^3\right)^\circ$, $Q_1(\zeta_0) = 2\leq Q_1(\zeta)$.

        \item[(2)] $\text{rank}(\zeta_0) = 2$.
        Since rank is a lower semi-continuous function, the set $\{\zeta\in\left(Z\otimes\mathbb{R}^3\right)^\circ \mid \text{rank}(\zeta)> 1\}$ is open. Thus there exists $\delta_0 >0$ such that for all $\zeta$ with $|\zeta - \zeta_0|<\delta_0$, we have $\text{rank}(\zeta)\geq 2$.
        \begin{itemize}
            \item Case: $\text{rank}(\zeta) = 2$. From subsection 4.2, there exists two planes $L_\zeta$ and $L_{\zeta_0}$ in $\mathfrak{h}^{\mathbb{C}}$ inducing two decompositions of the root system $\Phi$:
            \begin{align*}
                \Phi &= \Phi_1 \sqcup \cdots \sqcup\Phi_{s_{\zeta}}\\
                \Phi &= \Phi_1' \sqcup \cdots \sqcup\Phi_{s_{\zeta_0}}'
            \end{align*}
            which yield
            \begin{align*}
                Q_1(\zeta_0) &= \sum_{k = 1}^{s_{\zeta_0}}\#\,\{u_k^1,u_k^2\}\\
                Q_1(\zeta) &= \sum_{j = 1}^{s_{\zeta}}\#\,\{v_j^1,v_j^2\}
            \end{align*}
            
            First we show that $s_{\zeta_0}\leq s_{\zeta}$ for $\zeta$ sufficiently closed to $\zeta_0$ :  For each $k = 1,\cdots,s_{\zeta_0}$, select any $\alpha_k\in\Phi_k$ and denote the unique number $m(k)\in\{1,\cdots,s_{\zeta}\}$ such that $\alpha_k\in\Phi'_{m(k)}$. We prove by contradiction that the map (depending on choice of $\alpha_k$'s) $m: k\mapsto m(k)$ is injective: suppose there exists distinct $(k_1,k_2)\in \{1,\cdots,s_{\zeta_0}\}\cross\{1,\cdots,s_{\zeta_0}\}$ with $m(k_1) = m(k_2)$, that is 
            $
            \begin{cases}
                (\ker\alpha_{k_1}) \cap \,L_{\zeta_0} \neq (\ker\alpha_{k_2}) \cap \, L_{\zeta_0},\\
                (\ker\alpha_{k_1}) \cap \, L_{\zeta} = (\ker\alpha_{k_2}) \cap \,L_{\zeta}.
            \end{cases}
            $. 
            As $\zeta \to \zeta_0$, it's obvious that $L_{\zeta} \to L_{\zeta_0}$, so $(\ker\alpha_{k_i}) \cap L_{\zeta} \to (\ker\alpha_{k_i})\cap L_{\zeta_0}$ for $i = 1,2$, leading to the contradiction $(\ker \alpha_{k_1}) \cap L_{\zeta_0} = (\ker \alpha_{k_2}) \cap \,L_{\zeta_0}$.
            Therefore $m: k \mapsto m(k)$ is injective and then $s_{\zeta_0} \leq s_{\zeta}$.

            \quad To establish the inequality for $Q_1$, we analyze the the ramification of the map $p'_{\zeta_0}$, see subsection 4.2. 
            For $\alpha\in\Phi_k$, let $l_k := \ker\alpha_k\,\cap L_{\zeta_0}$ and $l_{m(k),\zeta} := \ker \alpha_k\, \cap L_{\zeta}$.  
            We claim that if $[l_{k}]$ is not a ramification value of $p'_{\zeta_0}$, then for $\zeta$ sufficiently closed to $\zeta_0$, $[l_{m(k),\zeta}]$ is not a ramification value of $p'_{\zeta}$: 
            Otherwise, there would exists a sequence $\zeta(n)$ with $|\zeta(n)-\zeta_0|<\frac{1}{n}$, such that $[l_{m(k),\zeta(n)}]$ is ramification of $p'_{\zeta(n)}$. Denoted by $T(n)$ the subset of preimage of $[l_{m(k),\zeta(n)}]$ under $p'_{\zeta(n)}$ where $\dd p'_{\zeta(n)}$ vanishes. Since $\mathbb{P}^1$ is compact, there exists $\{t_n\in T(n)\}$ such that $t_n$ converge to an accumulation point $t_0$. 
            Then take the limit
            \begin{equation*}
                \lim_{n\to\infty} \dd p'_{\zeta(n)}(t_n) = \dd p'_{\zeta_0}(t_0) = 0
            \end{equation*}
            Since $p'_{\zeta(n)}(T(n)) = [l_{m(k),\zeta(n)}]$ and $[l_{m(k),\zeta}] \to [l_{k}]$ when $\zeta\to \zeta_0$, 
            \begin{equation*}
                p'_{\zeta_0}(t_0) = \lim_{n\to \infty}p'_{\zeta(n)}(t_n) = \lim_{n \to \infty} [l_{m(k),\zeta(n)}] = [l_k].
            \end{equation*}
            It follows that $[l_k]$ is a ramification value of $p'_{\zeta_0}$, a contradiction. 
            Therefore, we see for each $k$, whether $[l_k]$ is a ramification value or not, there always holds
            \begin{equation*}
                \#\{u_k^1,u_k^2\} \leq 
                \#\{v_{m(k)}^1,v_{m(k)}^2\}
            \end{equation*}
            and hence:
            \begin{equation*}
                Q_1(\zeta_0) = \sum_{k}^{s_{\zeta_0}}\#\,\{u_k^1,u_k^2\} \leq \sum_{j = 1}^{s_{\zeta}} \#\,\{v_j^1,v_j^2\} = Q_1(\zeta).
            \end{equation*}

        \item Case: $\text{rank}(\zeta) = 3$. Let $C_\zeta$ denote the projectivization of the image of period map $\Tilde{p_\zeta}$, which is a conic curve in $\mathbb{P}(\mathfrak{h}^{\mathbb{C}})$, see subsection 4.3. 
        Then as $\zeta \to \zeta_0$, $C_\zeta \to \mathbb{P}(L_{\zeta_0})\subset\mathbb{P}(\mathfrak{h}^{\mathbb{C}})$. 
        We claim that for every $\alpha\in\Phi$, there exists $\epsilon(\alpha)>0$ such that $\mathbb{P}(\ker\alpha)$ is not tangent to $C_{\zeta}$ whenever $|\zeta-\zeta_0|<\epsilon(\alpha)$:
        If not, that is $\exists\alpha_0\in\Phi$ such that $\forall\,\epsilon>0$, $\exists\,\zeta(\epsilon)$ with $|\zeta(\epsilon)- \zeta_0|<\epsilon$, the $\mathbb{P}(\ker\alpha_0)$ is tangent to $C_{\zeta(\epsilon)}$. Then as $\epsilon \to 0$, $\zeta(\epsilon) \to \zeta$, $C_{\zeta(\epsilon)} \to \mathbb{P}(L_{\zeta_0})$. Hence $\mathbb{P}(\ker\alpha_0)$ tangents to $\mathbb{P}(L_{\zeta_0})$. This implies $\mathbb{P}(L_{\zeta_0}) = \mathbb{P}(\ker\alpha_0)$. This contradicts to the fact that the image of period map is not contained in any $\ker \alpha$, see subsection 4.2.

        \quad Hence if we take $\epsilon < \min\{\epsilon(\alpha),\alpha\in\Phi\}$, and $\zeta$ with $|\zeta-\zeta_0|<\epsilon$, then for any $\alpha\in\Phi$, $\mathbb{P}(\ker\alpha)$ intersects $C_\zeta$ transverse. Therefore: 
        \begin{align*}
             Q_1(\zeta) & = \text{the number of intersection points of $C_\zeta$ with the set $\bigcup_{\alpha\in\Phi}\mathbb{P}(\ker\alpha)$}\\
             & \geq 2s_{\zeta_0} \geq   \sum_{k}^{s_{\zeta_0}}\#\,\{u_k^1,u_k^2\} = Q_1(\zeta_0).
        \end{align*}

        \end{itemize}
        
        \item[(3)] $\text{rank}(\zeta_0) = 3$. Since $\{\zeta\mid \text{rank}(\zeta) >2\}$ is open, $\exists\,\delta_3 >0 $ such that $\forall \zeta$ with $|\zeta - \zeta_0|<\delta_3$, $\text{rank}(\zeta) = 3$. As $\zeta \to \zeta_0$, $C_{\zeta} \to C_{\zeta_0}$ in $\mathbb{P}(\mathfrak{h}_{\mathbb{C}})$. The transverse points $\gamma_k$ of $C_{\zeta_0}$ (defined in subsection 4.3) perturbs to transverse points again. Hence, $Q_1(\zeta)\geq Q_1(\zeta_0)$.

    \end{itemize}
\end{proof}

\section{Proof of Corollary \ref{cor_3}}
Hitchin proved that for any hyperK\"ahler ALE with group $\mathbb{Z}_{2}$, its twistor space cannot be K\"ahlerian due to the existence of holomorphic curves fiberwised in twistor space \cite{Hit81}.
Hence we need only consider the case where $\Gamma \neq\mathbb{Z}_{2}$. 
If $\mathrm{Z}(X_{\zeta})$ is K\"ahlerian, then its canonical bundle $K_{\mathrm{Z}}$ would be negative established by Hitchin \cite[p.150]{Hit81}. 
Hence, the line bundle $K_{\mathrm{Z}}$ would be of negative degree when restricted on any compact holomorphic curves in $\mathrm{Z}(X_{\zeta})$. 
The following lemma that $K_{\mathrm{Z}}\simeq\pi^{*}\mathcal{O}_{\mathbb{P}^{1}}(-4)$ and the property of pullback, imply the line bundle $K_{\mathrm{Z}}$ is trivial when restricted to any fiber of $\pi$. 
Theorem \ref{tm_2} implies that there is always at least one fiber containing a rational curve. The existence of such holomorphic curves leads to the contradiction. Therefore, $\mathrm{Z}(X_{\zeta})$ cannot be K\"ahlerian.\qed
   
\begin{lemma}
    Denote $\mathrm{Z}(M)$ as the twistor space of a hyperK\"ahler manifold $\left(M,g,I,J,K\right)$. Recall the natural holomorphic projection $\pi: \mathrm{Z}(M)\to\mathbb{P}^{1}$. Then, the canonical bundle $K_{\mathrm{Z}(M)}$ of the complex threefold $\mathrm{Z}(M)$ is given by:
    $$
    K_{\mathrm{Z}(M)} \simeq \pi^{*}\mathcal{O}_{\mathbb{P}^{1}}(-4).
    $$
\end{lemma}
\begin{proof}
    Denote the three K\"ahler forms of $\left(M,g,I,J,K\right)$ by $\left(\omega_{I},\omega_{J},\omega_{K}\right)$. The twisted vertical two-form, $\Omega \in H^{0}(\mathrm{Z}(X_{\zeta}),\,\Lambda^{2}T_{F}^{*}\otimes\pi^{*}\mathcal{O}_{\mathbb{P}^{1}}(2))$, as defined in \cite{HKLM}, is given by 
    \begin{align*}
        \Omega(u) &= (\omega_{I} + \mathrm{i}\omega_{J}) + u(2\omega_{I}) - u^{2}(\omega_{J} - \mathrm{i}\omega_{K}),\\
        \Omega(v) &= v^{2}(\omega_{I} + \mathrm{i}\omega_{J}) + v(2\omega_{I}) - (\omega_{J} - \mathrm{i}\omega_{K}),
    \end{align*}
    where $u \in \mathbb{C} \simeq U_{0}\subset \mathbb{P}^{1}$ and $v \in \mathbb{C} \simeq U_{1}\subset \mathbb{P}^{1}$. We now take a nowhere vanishing section of $K_{\mathrm{Z}(M)}$ over $U_{0}$, given by $\dd u \wedge \Omega(u)$, and a nowhere vanishing section of $K_{\mathrm{Z}(M)}$ over $U_{1}$, given by $\dd v \wedge \Omega(v)$. The transformation over $U_{0}\cap U_{1}$ is given by $\dd v\wedge\Omega(v) = -\frac{1}{u^{4}} \dd u \wedge \Omega(u) $. Therefore, by examining the transition function, we conclude that $ K_{\mathrm{Z}(M)}\simeq \pi^{*}\mathcal{O}_{\mathbb{P}^{1}}(-4).$
\end{proof}

\begin{example}
\label{ex_twi_Z2} 
For $\Gamma = \mathbb{Z}_{2}$, by Corollary \ref{cor_Z2}, the hyperK\"ahler ALE structure is realized on the minimal resolution $\widetilde{\mathbb{C}^{2}/\mathbb{Z}_{2}}$. As complex manifolds, there is an isomorphism 
\begin{equation*}
\widetilde{\mathbb{C}^{2}/\mathbb{Z}_{2}}\simeq\mathcal{O}_{\mathbb{P}^{1}}(-2)\simeq\mathcal{T}^{*}\,\mathbb{P}^{1}.    
\end{equation*}
The hyperK\"ahler metric on $\mathcal{O}_{\mathbb{P}^{1}}(-2)$ was constructed by Eguchi-Hanson \cite{EguchiHanson:78} and also by Calabi \cite{Calabi:79}. 
Its twistor space $\mathrm{Z}$, as explained in \cite{Hit79}, is obtained by resolving the singularities on a singular hypersurface $\Tilde{\mathrm{Z}}$ of a four-dimensional complex manifold:
\begin{equation*}
    \Tilde{\mathrm{Z}} := \left\{xy = (z-p_{1})(z-p_{2})\right\} \subset \mathcal{O}_{\mathbb{P}^1}(2)\oplus\mathcal{O}_{\mathbb{P}^1}(2)\oplus\mathcal{O}_{\mathbb{P}^1}(2)
\end{equation*}
where $x,y,z$ are parameters of the fibers of three line bundle, and $p_{1},p_{2}$ are some two distinct sections of $\mathcal{O}_{\mathbb{P}^1}(2)$.

\end{example}

$\mathbf{Acknowledgements.}$ 
Both authors express their profound gratitude to Friedrich Knop for his insightful responses to our inquiries on MathOverflow and for generously sharing the invaluable reference \cite{slo80}. The first author is particularly indebted to Song Sun for the enlightening discussions that substantially enriched this research.

\bibliographystyle{plain}
\bibliography{ref}

\begin{thebibliography}{10}

\bibitem{BHPV_book}
Wolf~P. Barth, Klaus Hulek, Chris A.~M. Peters, and Antonius Van~de Ven.
\newblock {\em Compact complex surfaces}, volume~4 of {\em Ergebnisse der Mathematik und ihrer Grenzgebiete. 3. Folge. A Series of Modern Surveys in Mathematics [Results in Mathematics and Related Areas. 3rd Series. A Series of Modern Surveys in Mathematics]}.
\newblock Springer-Verlag, Berlin, second edition, 2004.

\bibitem{Lie_Bour_46}
Nicolas Bourbaki.
\newblock {\em Lie groups and {L}ie algebras. {C}hapters 4--6}.
\newblock Elements of Mathematics (Berlin). Springer-Verlag, Berlin, 2002.
\newblock Translated from the 1968 French original by Andrew Pressley.

\bibitem{Calabi:79}
E.~Calabi.
\newblock {M{\'e}triques k{\"a}hl{\'e}riennes et fibr{\'e}s holomorphes}.
\newblock {\em Ann. Sci. {\'E}cole Norm. Sup. (4)}, 12(2):269--294, 1979.

\bibitem{EguchiHanson:78}
Tohru Eguchi and Adrew~J. Hanson.
\newblock {Asymptotically flat self--dual solutions to Euclidean gravity}.
\newblock {\em Phys. Lett. B}, 74:249--251, 1978.

\bibitem{HR84}
Hans Grauert and Reinhold Remmert.
\newblock {\em Coherent analytic sheaves}, volume 265 of {\em Grundlehren der mathematischen Wissenschaften [Fundamental Principles of Mathematical Sciences]}.
\newblock Springer-Verlag, Berlin, 1984.

\bibitem{Hit79}
N.~J. Hitchin.
\newblock Polygons and gravitons.
\newblock {\em Math. Proc. Cambridge Philos. Soc.}, 85(3):465--476, 1979.

\bibitem{Hit81}
N.~J. Hitchin.
\newblock K\"{a}hlerian twistor spaces.
\newblock {\em Proc. London Math. Soc. (3)}, 43(1):133--150, 1981.

\bibitem{HKLM}
N.~J. Hitchin, A.~Karlhede, U.~Lindstr\"{o}m, and M.~Ro\v{c}ek.
\newblock Hyper-{K}\"{a}hler metrics and supersymmetry.
\newblock {\em Comm. Math. Phys.}, 108(4):535--589, 1987.

\bibitem{Kas-Schlessinger_versal-deformation_1972}
Arnold Kas and Michael Schlessinger.
\newblock On the versal deformation of a complex space with an isolated singularity.
\newblock {\em Math. Ann.}, 196:23--29, 1972.

\bibitem{Kron1}
P.~B. Kronheimer.
\newblock The construction of {ALE} spaces as hyper-{K}\"{a}hler quotients.
\newblock {\em J. Differential Geom.}, 29(3):665--683, 1989.

\bibitem{Kron2}
P.~B. Kronheimer.
\newblock A {T}orelli-type theorem for gravitational instantons.
\newblock {\em J. Differential Geom.}, 29(3):685--697, 1989.

\bibitem{slodowy_LNM}
P.~Slodowy.
\newblock Platonic solids, {K}leinian singularities, and {L}ie groups.
\newblock In {\em Algebraic geometry ({A}nn {A}rbor, {M}ich., 1981)}, volume 1008 of {\em Lecture Notes in Math.}, pages 102--138. Springer, Berlin, 1983.

\bibitem{slo80}
Peter Slodowy.
\newblock {\em Simple singularities and simple algebraic groups}, volume 815 of {\em Lecture Notes in Mathematics}.
\newblock Springer, Berlin, 1980.

\bibitem{Tjurina_locally-semi-universal}
G.~N. Tjurina.
\newblock Locally semi-universal flat deformations of isolated singularities of complex spaces.
\newblock {\em Izv. Akad. Nauk SSSR Ser. Mat.}, 33:1026--1058, 1969.

\bibitem{Tyu70}
G.~N. Tjurina.
\newblock Resolution of singularities of flat deformations of double rational points.
\newblock {\em Funkcional. Anal. i Prilo\v zen.}, 4(1):77--83, 1970.

\end{thebibliography}

\end{document}